# $L^1$ HARDY INEQUALITIES WITH WEIGHTS*

GEORGIOS PSARADAKIS

psaradakis.georgios@gmail.com

**Abstract**

We prove sharp homogeneous improvements to $L^1$ weighted Hardy inequalities involving distance from the boundary. In the case of a smooth domain we obtain lower and upper estimates for the best constant of the remainder term. These estimates are sharp in the sense that they coincide when the domain is a ball or an infinite strip. In the case of a ball we also obtain further improvements.

# Contents

# 1 Introduction

Hardy's inequality involving distance from the boundary of a convex set $\Omega \subsetneq \mathbb{R}^n$; $n \geq 1$, asserts that

$$\int_\Omega |\nabla u|^p \mathrm{d}x \geq \Big(\frac{p-1}{p}\Big)^p \int_\Omega \frac{|u|^p}{d^p} \mathrm{d}x, \quad p > 1, \tag{1.1}$$

for all $u \in C_c^\infty(\Omega)$, where $d \equiv d(x) := \mathrm{dist}(x, \mathbb{R}^n \setminus \Omega)$. Due to [HLP], [D], [MS] and [MMP] the constant appearing in (1.1) is optimal. After the pioneering results in [Mz] and [BrM], a sequence of papers have improved (1.1) by adding extra terms on its right hand side, see for instance [BFT2], [BFT3], [FMT3], [FTT] and primarily [BFT1] and [FMT1], [FMT2] where it was also noted that (1.1) remains valid with the sharp constant in more general sets than convex ones, and in particular in sets that satisfy $-\Delta d \geq 0$ in the distributional sense.

In the case $p = 1$, (1.1) reduces to a trivial inequality, at least for sets having non positive distributional Laplacian of the distance function. However, in the one dimensional case, the following $L^1$ weighted Hardy inequality is well known:

$$\int_0^\infty \frac{|u'(x)|}{x^{s-1}} \mathrm{d}x \geq (s-1) \int_0^\infty \frac{|u(x)|}{x^s} \mathrm{d}x, \quad s > 1, \tag{1.2}$$

for all absolutely continuous functions $u : [0, \infty) \to \mathbb{R}$, such that $u(0) = 0$. This is the special case $p = 1$ of Theorem 330 in [HLP]. Inequality (1.2) is, in fact, an equality for $u$ increasing, and thus the constant on the right hand side is sharp.

In this work we are concerned with the higher-dimensional generalizations of (1.2). Let $\Omega \subsetneq \mathbb{R}^n (n \geq 2)$ be open and let $d \equiv d(x) := \mathrm{dist}(x, \mathbb{R}^n \setminus \Omega)$. We deal with inequalities of the type

$$\int_\Omega \frac{|\nabla u|}{d^{s-1}} \mathrm{d}x \geq \mathcal{B}_0 \int_\Omega \frac{|u|}{d^s} \mathrm{d}x + \mathcal{B} \int_\Omega V(d)|u| \mathrm{d}x, \quad s \geq 1, \tag{1.3}$$

valid for all $u \in C_c^\infty(\Omega)$. Here $V$ is a potential function, i.e., nonnegative and $V \in L^1_{loc}(\mathbb{R}^+)$, and $\mathcal{B}_0 \geq 0$, $\mathcal{B} \in \mathbb{R}$. Questions concerning sets for which this inequality is valid, sharp constants, possible improvements and optimal potentials will be studied. Our first Theorem reads as follows

**Theorem A** *Let $\Omega$ be a domain in $\mathbb{R}^n$ with boundary of class $\mathcal{C}^2$ satisfying a uniform interior sphere condition, and we denote by $\underline{\mathcal{H}}$ the infimum of the mean curvature of the boundary. Then there exists $\mathcal{B}_1 \geq (n-1)\underline{\mathcal{H}}$ such that for all $u \in C_c^\infty(\Omega)$ and all $s \geq 1$*

$$\int_\Omega \frac{|\nabla u|}{d^{s-1}} \mathrm{d}x \geq (s-1) \int_\Omega \frac{|u|}{d^s} \mathrm{d}x + \mathcal{B}_1 \int_\Omega \frac{|u|}{d^{s-1}} \mathrm{d}x. \tag{1.4}$$

*Let $s \geq 2$. If $\Omega$ is a bounded domain in $\mathbb{R}^n$ with boundary of class $\mathcal{C}^2$ having strictly positive mean curvature, then the constant $s - 1$ in the first term as well as the exponent $s - 1$ on the*

*distance function on the remainder term in (1.4), are optimal. In addition, we have the following estimates*

$$(n-1)\underline{\mathcal{H}} \leq \mathcal{B}_1 \leq \frac{n-1}{|\partial\Omega|}\int_{\partial\Omega}\mathcal{H}(y)\mathrm{d}S_y, \tag{1.5}$$

*where $\mathcal{H}(y)$ is the mean curvature of the boundary at $y \in \partial\Omega$, and $\underline{\mathcal{H}}$ is its minimum value.*

The following result, which is of independent interest, played a key role in establishing Theorem A

**Theorem B** *Let $\Omega \subset \mathbb{R}^n$ be a domain with boundary of class $\mathcal{C}^2$ satisfying a uniform interior sphere condition. Then $\mu := (-\Delta d)\mathrm{d}x$ is a signed Radon measure on $\Omega$. Let $\mu = \mu_{ac} + \mu_s$ be the Lebesgue decomposition of $\mu$ with respect to $\mathcal{L}^n$, i.e. $\mu_{ac} \ll \mathcal{L}^n$ and $\mu_s \perp \mathcal{L}^n$. Then $\mu_s \geq 0$ in $\Omega$, and $\mu_{ac} \geq (n-1)\underline{\mathcal{H}}\mathrm{d}x$ a.e. in $\Omega$, where $\underline{\mathcal{H}} := \inf_{y \in \partial\Omega}\mathcal{H}(y)$.*

For domains with boundary of class $\mathcal{C}^2$ satisfying a uniform interior sphere condition, $-\Delta d$ is a continuous function in a tubular neighborhood of the boundary and, moreover, $-\Delta d(y) = (n-1)\mathcal{H}(y)$ for any $y \in \partial\Omega$. This fact together with Theorem B leads to

**Corollary** *Let $\Omega$ be a domain with boundary of class $\mathcal{C}^2$ satisfying a uniform interior sphere condition. Then $\Omega$ is mean convex, i.e., $\mathcal{H}(y) \geq 0$ for all $y \in \partial\Omega$, if and only if $-\Delta d \geq 0$ holds in $\Omega$, in the sense of distributions.*

We note that a set $\Omega \subsetneq \mathbb{R}^n$ with distance function having non positive distributional Laplacian, is shown in [4-5] and [13-15] to be the natural geometric assumption for the validity of various Hardy inequalities.

In special geometries we are able to compute the best constant $\mathcal{B}_1$ in (1.4):

In case $\Omega$ is a ball of radius $R$ then the upper and lower estimates (1.5) coincide, yielding $\mathcal{B}_1 = (n-1)/R$. One then may ask whether (1.4) can be further improved. We provide a full answer to this question by showing that for $s \geq 2$ one can add a finite series of $[s] - 1$ terms on the right hand side before adding an optimal logarithmic correction. More precisely we prove the following

**Theorem C** *Let $B_R$ be a ball of radius $R$. Then, $(i)$ For all $u \in C_c^\infty(B_R)$, all $s \geq 2$, $\gamma > 1$, it holds that*

$$\int_{B_R}\frac{|\nabla u|}{d^{s-1}}\mathrm{d}x \geq (s-1)\int_{B_R}\frac{|u|}{d^s}\mathrm{d}x + \sum_{k=1}^{[s]-1}\frac{n-1}{R^k}\int_{B_R}\frac{|u|}{d^{s-k}}\mathrm{d}x + \frac{C}{R^{s-1}}\int_{B_R}\frac{|u|}{d}X^\gamma\Big(\frac{d}{R}\Big)\mathrm{d}x, \tag{1.6}$$

*where $X(t) := (1 - \log t)^{-1}$, $t \in (0,1]$ and $C \geq \gamma - 1$. The exponents $s$ and $s - k$; $k = 1, 2, ..., [s] - 1$, on the distance function, as well as the constants $s - 1$, $(n-1)/R^k$; $k = 1, 2, ..., [s] - 1$, in the first and the summation terms, respectively, are optimal. The last term in (1.6) is optimal in the sense that if $\gamma = 1$, there is not positive constant $C$ such that (1.6) holds.*

*(ii) For all $u \in C_c^\infty(B_R)$, all $1 \leq s < 2$, $\gamma > 1$, it holds that*

$$\int_{B_R} \frac{|\nabla u|}{d^{s-1}} dx \geq (s-1) \int_{B_R} \frac{|u|}{d^s} dx + \frac{C}{R^{s-1}} \int_{B_R} \frac{|u|}{d} X^\gamma\left(\frac{d}{R}\right) dx, \tag{1.7}$$

*where $X(t) := (1 - \log t)^{-1}$, $t \in (0, 1]$ and $C \geq \gamma - 1$. The last term in (1.7) is optimal in the sense that if $\gamma = 1$, there is not positive constant $C$ such that (1.7) holds.*

Note that this is in contrast with the results in case $p > 1$, where an infinite series involving optimal logarithmic terms can be added (see [BFT2]) and ([BFT3]).

In case $\Omega$ is an infinite strip, using a more general upper bound on $\mathcal{B}_1$ (see Theorem 4.7), we prove that $\mathcal{B}_1 = 0$. As a matter of fact the finite series structure of (1.6) disappears and only the final logarithmic correction term survives. More precisely

**Theorem D** *Let $S_R$ be an infinite strip of inner radius $R$. For all $u \in C_c^\infty(S_R)$, all $s \geq 1$, $\gamma > 1$, there holds*

$$\int_{S_R} \frac{|\nabla u|}{d^{s-1}} dx \geq (s-1) \int_{S_R} \frac{|u|}{d^s} dx + \frac{C}{R^{s-1}} \int_{S_R} \frac{|u|}{d} X^\gamma\left(\frac{d}{R}\right) dx, \tag{1.8}$$

*where $C \geq \gamma - 1$. The last term in (1.8) is optimal in the sense that if $\gamma = 1$, there is not positive constant $C$ such that (1.8) holds.*

The paper is organized as follows. In §2, recalling the semiconcavity properties of the distance function, we prove weighted $L^1$ Hardy inequalities in sets without regularity assumptions on the boundary. General open sets, sets with non negative distributional Laplacian of the distance function, as well as sets with positive reach are considered. Remainders for sets having finite inner radius are obtained in the first two cases and extremal domains are given. The results imply in particular inequality (1.8). In §3, after recalling further properties of the distance function for smooth domains, we prove Theorem B. Theorem A and the optimality in Theorem D is then proved in §4, where also an interesting lower bound for the Cheeger constant of smooth, strictly mean convex domains is deduced (see Corollary 4.6). In §5, Theorem C is proved and in the final section we discuss $L^p$ analogs of our results.

After this work was completed we found that Corollary following Theorem B of this introduction is also noted in [LL]. It turns out this is originally due to Gromov (see [Gr]-pg 18-19). For proofs of this corollary (different from the one in this paper) see [LLL] and [Gr]-§5.

## 2 Inequalities in sets without regularity assumptions on the boundary

Since all inequalities of this paper will follow by the integration by parts formula, we formalize it as follows: let $\Omega$ be an open set in $\mathbb{R}^n$ and $T$ be a vector field on $\Omega$. Integrating by parts and using elementary inequalities, we get

$$\int_\Omega |T||\nabla u|\mathrm{d}x \geq \int_\Omega \mathrm{div}(T)|u|\mathrm{d}x, \tag{2.1}$$

for all $u \in C_c^\infty(\Omega)$, where we have also used the fact that $|\nabla|u|| = |\nabla u|$ a.e. in $\Omega$.

## 2.1 General sets

In this subsection we recall some properties of the distance function to the boundary of a general open set and then prove various weighted $L^1$ Hardy inequalities.

Let $\Omega \subsetneq \mathbb{R}^n$ be open. We set $d : \mathbb{R}^n \to [0, \infty)$ by $d(x) := \inf\{|x - y| : y \in \mathbb{R}^n \setminus \Omega\}$. It is well known that $d$ is Lipschitz continuous on $\mathbb{R}^n$ and in particular $|\nabla d(x)| = 1$ a.e. in $\Omega$ (see [F]-Theorem 4.8). The next property of $d$ can be found for example in [CS]-Proposition 2.2.2.(i) & Proposition 1.1.3.(c),(e). We prove it for completeness.

**Lemma 2.1.** *Let $\Omega \subsetneq \mathbb{R}^n$ be open. It holds that*

$$-d\Delta d \geq -(n-1) \text{ in } \Omega \text{ in the sense of distributions.} \tag{2.2}$$

*Proof.* Estimate (2.2) rests on the fact that the function $A : \mathbb{R}^n \to \mathbb{R}$ defined by $A(x) := |x|^2 - d^2(x)$ is convex (see also [AmbM]-§2). To see this, we take $x \in \mathbb{R}^n$ and let $y \in \mathbb{R}^n$ be such that $d(x) = |x - y|$. For any $z \in \mathbb{R}^n$ we get

$$\begin{aligned}
A(x+z) + A(x-z) - 2A(x) &= 2|z|^2 - (d^2(x-z) + d^2(x+z) - 2d^2(x)) \\
&\geq 2|z|^2 - (|x+z-y|^2 + |x-z-y|^2 - 2|x-y|^2) \\
&= 0.
\end{aligned}$$

Since $A(x)$ is also continuous, we obtain that $A(x)$ is convex (see [CS]-Proposition A1.2). It follows by [EvG]-§6.3-Theorem 2, that the distributional Laplacian of $A$ is a nonnegative Radon measure on $\mathbb{R}^n$. Since in $\Omega$ we have $\Delta A = 2(n - 1 - d\Delta d)$ in the sense of distributions, the result follows. ∎

The weighted $L^1$ Hardy inequalities we obtain are deduced from the following basic fact

**Lemma 2.2.** *Let $\Omega \subsetneq \mathbb{R}^n$ be open. For all $u \in C_c^\infty(\Omega)$ and all $s \geq 1$*

$$\int_\Omega \frac{|\nabla u|}{d^{s-1}}\mathrm{d}x \geq (s-1)\int_\Omega \frac{|u|}{d^s}\mathrm{d}x + \int_\Omega \frac{|u|}{d^{s-1}}(-\Delta d)\mathrm{d}x, \tag{2.3}$$

*where $-\Delta d$ is meant in the distributional sense. If $\Omega$ is bounded, then equality holds for $u_\varepsilon(x) = (d(x))^{s-1+\varepsilon} \in W_0^{1,1}(\Omega; d^{-(s-1)})$, $\varepsilon > 0$.*

*Proof.* Inequality (2.3) follows from (2.1) by setting $T(x) = -(d(x))^{1-s}\nabla d(x)$ for a.e. $x \in \Omega$, while the second statement is easily checked. ∎

A covering of $\Omega$ by cubes was used in [Avkh] to prove the next Theorem. We present an elementary proof.

**Theorem 2.3.** *Let $\Omega \subsetneq \mathbb{R}^n$ be open. For all $u \in C_c^\infty(\Omega)$ and all $s > n$, it holds that*

$$\int_\Omega \frac{|\nabla u|}{d^{s-1}} \mathrm{d}x \geq (s-n) \int_\Omega \frac{|u|}{d^s} \mathrm{d}x. \tag{2.4}$$

*Proof.* Coupling (2.2) and (2.3), we get

$$\begin{aligned}
\int_\Omega \frac{|\nabla u|}{d^{s-1}} \mathrm{d}x &\geq (s-1) \int_\Omega \frac{|u|}{d^s} \mathrm{d}x - (n-1) \int_\Omega \frac{|u|}{d^s} \mathrm{d}x \\
&= (s-n) \int_\Omega \frac{|u|}{d^s} \mathrm{d}x.
\end{aligned}$$
∎

**Remark 2.4.** The constant appearing on the right hand side of (2.4) is just a lower bound for the best constant. The best constant in (2.4) differs from one open set to another. However, $\mathbb{R}^n \setminus \{0\}$ serves as an extremal domain for Theorem 2.3. More precisely, letting $\Omega = \mathbb{R}^n \setminus \{0\}$, we have $d(x) = |x|$ and (2.4) reads as follows

$$\int_{\mathbb{R}^n} \frac{|\nabla u|}{|x|^{s-1}} \mathrm{d}x \geq (s-n) \int_{\mathbb{R}^n} \frac{|u|}{|x|^s} \mathrm{d}x, \quad s > n, \tag{2.5}$$

for all $u \in C_c^\infty(\mathbb{R}^n \setminus \{0\})$. To illustrate the optimality of the constant on the right hand side of (2.5), we define the following function

$$u_\delta(x) := \chi_{B_\eta \setminus B_\delta}(x), \quad x \in \mathbb{R}^n, \tag{2.6}$$

where, for any $r > 0$, by $B_r$ we henceforth denote the open ball of radius $r$ with center at the origin. Here $0 < \delta < \eta$ and $\eta$ is fixed. The distributional gradient of $u_\delta$ is $\nabla u_\delta = \vec{\nu}_{\partial B_\delta} \delta_{\partial B_\delta} - \vec{\nu}_{\partial B_\eta} \delta_{\partial B_\eta}$ where, for any $r > 0$, $\vec{\nu}_{\partial B_r}$ stands for the outward pointing unit normal vector field along $\partial B_r = \{x \in \mathbb{R}^n : |x| = r\}$, and by $\delta_{\partial B_r}$ we denote the Dirac measure on $\partial B_r$. Moreover, the total variation of $\nabla u_\delta$ is $|\nabla u_\delta| = \delta_{\partial B_\delta} + \delta_{\partial B_\eta}$. Using the co-area formula, we get

$$\begin{aligned}
\frac{\int_{\mathbb{R}^n} \frac{|\nabla u_\delta|}{|x|^{s-1}} \mathrm{d}x}{\int_{\mathbb{R}^n} \frac{|u_\delta|}{|x|^s} \mathrm{d}x} &= \frac{\delta^{1-s} |\partial B_\delta| + \eta^{1-s} |\partial B_\eta|}{\int_\delta^\eta r^{-s} |\partial B_r| \mathrm{d}r} \\
&= \frac{\delta^{n-s} + \eta^{n-s}}{\int_\delta^\eta r^{n-s-1} \mathrm{d}r} \\
&= (s-n) \frac{\delta^{n-s} + \eta^{n-s}}{\delta^{n-s} - \eta^{n-s}} \\
&\to s-n, \quad \text{as } \delta \downarrow 0.
\end{aligned}$$

Although not smooth, functions like $u_\delta$ defined in (2.6) belong to $BV(\mathbb{R}^n)$ (the space of functions of bounded variation in $\mathbb{R}^n$), and thus we can use a $C_c^\infty$ approximation so that the calculation above to hold in the limit (see for instance [EvG]-§5.2).

**Theorem 2.5.** *Let $\Omega \subsetneq \mathbb{R}^n$ be open and such that $R := \sup_{x \in \Omega} d(x) < \infty$. For all $u \in C_c^\infty(\Omega)$, all $s \geq n, \gamma > 1$, it holds that*

$$\int_\Omega \frac{|\nabla u|}{d^{s-1}} \mathrm{d}x \geq (s-n) \int_\Omega \frac{|u|}{d^s} \mathrm{d}x + \frac{C}{R^{s-n}} \int_\Omega \frac{|u|}{d^n} X^\gamma\left(\frac{d}{R}\right) \mathrm{d}x, \tag{2.7}$$

*where $C \geq \gamma - 1$.*

*Proof.* We set $T(x) = -(d(x))^{1-s}[1 - (d(x)/R)^{s-n} X^{\gamma-1}(d(x)/R)] \nabla d(x)$ for a.e. $x \in \Omega$. Since $|1 - (d(x)/R)^{s-n} X^{\gamma-1}(d(x)/R)| \leq 1$ for all $x \in \Omega$, we have

$$\int_\Omega |T||\nabla u| \mathrm{d}x \leq \int_\Omega \frac{|\nabla u|}{d^{s-1}} \mathrm{d}x.$$

Using the rule $\nabla X^{\gamma-1}(d(x)/R) = (\gamma-1) X^\gamma(d(x)/R) \frac{\nabla d(x)}{d(x)}$ for a.e. $x \in \Omega$, we compute

$$\begin{aligned}
\mathrm{div}(T) &= (s-1) d^{-s}[1 - (d/R)^{s-n} X^{\gamma-1}(d/R)] + \frac{s-n}{R^{s-n}} d^{-n} X^{\gamma-1}(d/R) \\
&\quad + \frac{\gamma-1}{R^{s-n}} d^{-n} X^\gamma(d/R) + d^{1-s}[1 - (d/R)^{s-n} X^{\gamma-1}(d/R)](-\Delta d).
\end{aligned}$$

Since $1 - (d(x)/R)^{s-n} X^{\gamma-1}(d(x)/R) \geq 0$ for all $x \in \Omega$, we use (2.2) on the last term of the above equality and a straightforward computation gives

$$\mathrm{div}(T) \geq (s-n) d^{-s} + \frac{\gamma-1}{R^{s-n}} d^{-n} X^\gamma(d/R).$$

This means that

$$\int_\Omega \mathrm{div}(T) |u| \mathrm{d}x \geq (s-n) \int_\Omega \frac{|u|}{d^s} \mathrm{d}x + \frac{\gamma-1}{R^{s-n}} \int_\Omega \frac{|u|}{d^n} X^\gamma(d/R) \mathrm{d}x.$$

and the result follows from (2.1). ∎

**Remark 2.6.** A punctured domain serves as an extremal domain for Theorem 2.5. More precisely, let $\Omega = U \setminus \{0\}$ where $U$ is an open, connected subset of $\mathbb{R}^n$ containing the origin and satisfying $R := \sup_{x \in U} d(x) < \infty$. We define $u_\delta$ as in (2.6), where $\eta$ is fixed and sufficiently small such that $d(x) = |x|$ in $B_\eta$. For any $s \geq n$, we have

$$\begin{aligned}
\frac{\int_\Omega \frac{|\nabla u_\delta|}{|x|^{s-1}} \mathrm{d}x - (s-n) \int_\Omega \frac{|u_\delta|}{|x|^s} \mathrm{d}x}{\int_\Omega \frac{|u_\delta|}{|x|^n} X(|x|/R) \mathrm{d}x} &= \frac{\delta^{1-s}|\partial B_\delta| + \eta^{1-s}|\partial B_\eta| - (s-n) \int_\delta^\eta r^{-s} |\partial B_r| \mathrm{d}r}{\int_\delta^\eta r^{-n} X(r/R) |\partial B_r| \mathrm{d}r} \\
&= \frac{\delta^{n-s} + \eta^{n-s} - (s-n) \int_\delta^\eta r^{n-s-1} \mathrm{d}r}{\int_\delta^\eta r^{-1} X(r/R) \mathrm{d}r} \\
&= \frac{2\eta^{n-s}}{\log\left(\frac{X(\eta/R)}{X(\delta/R)}\right)} \\
&= o_\delta(1).
\end{aligned}$$

Thus, for a punctured domain inequality (2.7) does not hold when $\gamma = 1$, as well as the exponent $n$ on the second term of the right hand side in (2.7) cannot be increased.

**Theorem 2.7.** *Let $\Omega \subsetneq \mathbb{R}^n$ be open and such that $R := \sup_{x \in \Omega} d(x) < \infty$. For all $u \in C_c^\infty(\Omega)$ and all $s > n$, it holds that*

$$\int_\Omega \frac{|\nabla u|}{d^{s-1}} \mathrm{d}x - (s-n) \int_\Omega \frac{|u|}{d^s} \mathrm{d}x \geq \frac{1}{R^{s-n}} \int_\Omega \frac{|\nabla u|}{d^{n-1}} \mathrm{d}x. \tag{2.8}$$

**Proof.** We set $\vec{T}(x) = -(d(x))^{1-s}[1 - (d(x)/R)^{s-n}] \nabla d(x)$ for a.e. $x \in \Omega$. Since

$$|\vec{T}(x)| = (d(x))^{1-s}\left[1 - \left(\frac{d(x)}{R}\right)^{s-n}\right] \quad \text{a.e. } x \in \Omega,$$

we have

$$\int_\Omega |\vec{T}||\nabla u| \mathrm{d}x = \int_\Omega \frac{|\nabla u|}{d^{s-1}} \mathrm{d}x - \frac{1}{R^{s-n}} \int_\Omega \frac{|\nabla u|}{d^{n-1}} \mathrm{d}x.$$

We also calculate

$$\operatorname{div}(\vec{T}) = (s-1)d^{-s}[1-(d/R)^{s-n}] + \frac{s-n}{R^{s-n}}d^{-n} + d^{1-s}[1-(d/R)^{s-n}](-\Delta d), \quad \text{in } \Omega,$$

in the distributional sense. Since $1 - (d(x)/R)^{s-n} \geq 0$ for all $x \in \Omega$, we may use Lemma 2.1 on the last term of the above equality and after a straightforward computation to obtain

$$\int_\Omega \operatorname{div}(\vec{T})|u| \mathrm{d}x \geq (s-n) \int_\Omega \frac{|u|}{d^s} \mathrm{d}x.$$

The result follows from (2.1). ∎

**Remark 2.8.** A punctured domain serves also as an extremal domain for Theorem 2.7. As before, letting $\Omega = U \setminus \{0\}$, where $U$ is an open, connected subset of $\mathbb{R}^n$ containing the origin and satisfying $R := \sup_{x \in U} d(x) < \infty$, we define $u_\delta$ as in (2.6) where $\eta$ is fixed and sufficiently small such that $d(x) = |x|$ in $B_\eta$. By the co-area formula, for any $\varepsilon \geq 0$ we have

$$\begin{aligned}
\frac{\int_\Omega \frac{|\nabla u_\delta|}{|x|^{s-1}} \mathrm{d}x - (s-n) \int_\Omega \frac{|u_\delta|}{|x|^s} \mathrm{d}x}{\int_\Omega \frac{|\nabla u_\delta|}{|x|^{n-1+\varepsilon}} \mathrm{d}x} &= \frac{\delta^{1-s}|\partial B_\delta| + \eta^{1-s}|\partial B_\eta| - (s-n) \int_\delta^\eta r^{-s}|\partial B_r| \mathrm{d}r}{\delta^{1-n-\varepsilon}|\partial B_\delta| + \eta^{1-n-\varepsilon}|\partial B_\eta|} \\
&= \frac{\delta^{n-s} + \eta^{n-s} + \int_\delta^\eta (r^{n-s})' \mathrm{d}r}{\delta^{-\varepsilon} + \eta^{-\varepsilon}} \\
&= \frac{2\eta^{n-s}}{\delta^{-\varepsilon} + \eta^{-\varepsilon}} \\
&= \begin{cases} o_\delta(1) & \text{if } \varepsilon > 0 \\ \eta^{n-s} & \text{if } \varepsilon = 0. \end{cases}
\end{aligned}$$

Note that if $\varepsilon = 0$ then $\eta^{n-s} \downarrow R^{n-s}$ as $\eta \uparrow R$. □

## 2.2 Sets with $-\Delta d \geq 0$ in the sense of distributions

In this subsection we assume that

$$-\Delta d \geq 0 \text{ in } \Omega, \text{ in the sense of distributions.} \tag{C}$$

This condition was first used in the context of Hardy inequalities in [3-4] and has been used intensively in [13-15]. As we will prove in §3, domains with sufficiently smooth boundary carrying condition (C) are characterized as domains with nonnegative mean curvature of their boundary. However, in this section we do not impose regularity on the boundary.

**Theorem 2.9.** *Let $\Omega \subsetneq \mathbb{R}^n$ be open and such that condition* (C) *holds. For all $u \in C_c^\infty(\Omega)$ and all $s > 1$, it holds that*

$$\int_\Omega \frac{|\nabla u|}{d^{s-1}} \mathrm{d}x \geq (s-1) \int_\Omega \frac{|u|}{d^s} \mathrm{d}x. \tag{2.9}$$

*Moreover, the constant appearing on the right hand side of* (2.9) *is sharp.*

*Proof.* Since (C) holds we may cancel the last term in (2.3) and (2.9) follows. To prove the sharpness of the constant, we pick $y \in \partial\Omega$ and define the family of $W_0^{1,1}(\Omega; d^{-(s-1)})$ functions by $u_\varepsilon(x) := \phi(x)(d(x))^{s-1+\varepsilon}$, $\varepsilon > 0$, where $\phi \in C_c^\infty(B_\delta(y))$, $0 \leq \phi \leq 1$ and $\phi \equiv 1$ in $B_{\delta/2}(y)$, for some small but fixed $\delta$. We have

$$\begin{aligned}
\frac{\int_\Omega \frac{|\nabla u_\varepsilon|}{d^{s-1}} \mathrm{d}x}{\int_\Omega \frac{|u_\varepsilon|}{d^s} \mathrm{d}x} &\leq s - 1 + \varepsilon + \frac{\int_\Omega |\nabla \phi| d^\varepsilon \mathrm{d}x}{\int_\Omega \phi d^{-1+\varepsilon} \mathrm{d}x} \\
&\leq s - 1 + \varepsilon + \frac{C}{\int_{\Omega \cap B_{\delta/2}(y)} d^{-1+\varepsilon} \mathrm{d}x}, \\
&\leq s - 1 + o_\varepsilon(1)
\end{aligned}$$

where $C$ is some universal constant (not depending on $\varepsilon$). ∎

**Remark 2.10.** In view of Theorem 2.9 and Lemma 2.2, we see that if $\Omega$ is bounded and condition (C) holds, then all constants appearing in (2.3) are optimal.

**Theorem 2.11.** *Let $\Omega \subsetneq \mathbb{R}^n$ be open and such that condition* (C) *holds. Suppose in addition that $R := \sup_{x \in \Omega} d(x) < \infty$. For all $u \in C_c^\infty(\Omega)$, all $s \geq 1$, $\gamma > 1$, it holds that*

$$\int_\Omega \frac{|\nabla u|}{d^{s-1}} \mathrm{d}x \geq (s-1) \int_\Omega \frac{|u|}{d^s} \mathrm{d}x + \frac{C}{R^{s-1}} \int_\Omega \frac{|u|}{d} X^\gamma\left(\frac{d}{R}\right) \mathrm{d}x, \tag{2.10}$$

*where $C \geq \gamma - 1$.*

*Proof.* We set $T(x) = -(d(x))^{1-s}[1 - (d(x)/R)^{s-1} X^{\gamma-1}(d(x)/R)] \nabla d(x)$ for a.e. $x \in \Omega$. Since $|1 - (d(x)/R)^{s-1} X^{\gamma-1}(d(x)/R)| \leq 1$ for all $x \in \Omega$, we have

$$\int_\Omega |T| |\nabla u| \mathrm{d}x \leq \int_\Omega \frac{|\nabla u|}{d^{s-1}} \mathrm{d}x.$$

Using the rule $\nabla X^{\gamma-1}(d(x)/R) = (\gamma-1)X^\gamma(d(x)/R)\frac{\nabla d(x)}{d(x)}$ for a.e. $x \in \Omega$, by a straightforward calculation we arrive at

$$\int_\Omega \operatorname{div}(T)|u|\mathrm{d}x = (s-1)\int_\Omega \frac{|u|}{d^s}\mathrm{d}x + \frac{\gamma-1}{R^{s-1}}\int_\Omega \frac{|u|}{d}X^\gamma(d/R)\mathrm{d}x$$
$$+ \int_\Omega \frac{|u|}{d^{s-1}}[1 - (d/R)^{s-1}X^{\gamma-1}(d/R)](-\Delta d)\mathrm{d}x.$$

Since $1 - (d(x)/R)^{s-1}X^{\gamma-1}(d(x)/R) \geq 0$ for all $x \in \Omega$ and also (C) holds, we may cancel the last term and the result follows by (2.1). ∎

**Remark 2.12.** We prove in §4-Example 4.10 that an infinite strip is an extremal domain for Theorem 2.11. More precisely, if $\Omega = \{x = (x', x_n) : x' \in \mathbb{R}^{n-1}, 0 < x_n < 2R\}$ for some $R > 0$, then (2.10) fails for $\gamma = 1$ and thus the exponent 1 on the distance to the boundary in the remainder term of (2.10) cannot be increased.

The counterpart of Theorem 2.7 reads as follows

**Theorem 2.13.** *Let $\Omega \subsetneq \mathbb{R}^n$ be open, satisfies condition (C) and $R := \sup_{x \in \Omega} d(x) < \infty$. Then for all $u \in C_c^\infty(\Omega)$ and all $s > 1$*

$$\int_\Omega \frac{|\nabla u|}{d^{s-1}}\mathrm{d}x \geq (s-1)\int_\Omega \frac{|u|}{d^s}\mathrm{d}x + \frac{1}{R^{s-1}}\int_\Omega |\nabla u|\mathrm{d}x. \tag{2.11}$$

**Proof.** We insert $\vec{T}(x) = -(d(x))^{1-s}[1 - (d(x)/R)^{s-1}]\nabla d(x)$; a.e. $x \in \Omega$, in (2.1). Since $(d(x)/R)^{s-1} \leq 1$ for all $x \in \Omega$ we have

$$\int_\Omega |\vec{T}||\nabla u|\mathrm{d}x = \int_\Omega \frac{|\nabla u|}{d^{s-1}}\mathrm{d}x - \frac{1}{R^{s-1}}\int_\Omega |\nabla u|\mathrm{d}x.$$

On the other hand

$$\int_\Omega \operatorname{div}(\vec{T})|u|\mathrm{d}x = (s-1)\int_\Omega \frac{|u|}{d^s}\mathrm{d}x + \int_\Omega \frac{|u|}{d^{s-1}}(1 - (d/R)^{s-1})(-\Delta d)\mathrm{d}x$$
$$\geq (s-1)\int_\Omega \frac{|u|}{d^s}\mathrm{d}x,$$

where now we have used again the fact that $(d(x)/R)^{s-1} \leq 1$ for all $x \in \Omega$ and also (C). The result follows. ∎

An infinite strip is an extremal domain for (2.11), in the following sense

**Lemma 2.14.** *For fixed $R > 0$, set $S := \{x = (x', x_n) : x' \in \mathbb{R}^{n-1}, 0 < x_n < 2R\}$. Suppose that for some nonnegative $\alpha$ and $s \geq 1$, there holds*

$$\mathcal{C} := \inf_{u \in C_c^\infty(S)\setminus\{0\}} \tilde{Q}[u] \geq C_0 > 0,$$

*where*

$$\tilde{Q}[u] := \frac{\int_S \frac{|\nabla u|}{d^{s-1}} dx - (s-1) \int_S \frac{|u|}{d^s} dx}{\int_S \frac{|\nabla u|}{d^\alpha} dx}.$$

*Then $\alpha = 0$.*

**Proof.** For $s = 1$ it is obvious. Note also that it is enough to assume that $0 < \alpha < 1$. Let $s > 1$. Pick any $\phi \equiv \phi(x') \in C_c^1(\mathbb{R}^{n-1})$ such that $\text{sprt}\{\phi\} \subset B_1$, where $B_1$ is the $n-1$ dimensional open ball with radius 1 centered at $0'$. Let $\delta > 0$ and set $\phi_\delta \equiv \phi_\delta(x') := \phi(\delta x')$. Let also $0 < \varepsilon < \eta \leq R$. We test $\mathcal{C}$ with $u_{\varepsilon,\delta}(x) := \chi_{(\varepsilon,\eta)}(x_n)\phi_\delta(x')$. First note that

$$\nabla u_{\varepsilon,\delta}(x) = (\chi_{(\varepsilon,\eta)}(x_n)\nabla_{x'}\phi_\delta(x'), (\delta(x_n - \varepsilon) - \delta(x_n - \eta))\phi_\delta(x')),$$

where $\nabla_{x'} = (\frac{\partial}{\partial x_1}, \frac{\partial}{\partial x_2}, ..., \frac{\partial}{\partial x_{n-1}})$. Thus

$$|\nabla u_{\varepsilon,\delta}(x)| = \chi_{(\varepsilon,\eta)}(x_n)|\nabla_{x'}\phi_\delta(x')| + (\delta(x_n - \varepsilon) + \delta(x_n - \eta))|\phi_\delta(x')|.$$

Since $\eta \leq R/2$ we may substitute $d(x)$ by $x_n$ in $\tilde{Q}[u]$, and so

$$\begin{aligned}
\tilde{Q}[u_{\varepsilon,\delta}] &= \frac{\int_S \frac{|\nabla u_{\varepsilon,\delta}|}{x_n^{s-1}} dx - (s-1) \int_S \frac{|u_{\varepsilon,\delta}|}{x_n^s} dx}{\int_S \frac{|\nabla u_{\varepsilon,\delta}|}{x_n^\alpha} dx} \\
&= \frac{\int_\varepsilon^\eta \int_{B_{1/\delta}} \frac{|\nabla_{x'}\phi_\delta|}{x_n^{s-1}} dx' dx_n + (\frac{1}{\varepsilon^{s-1}} + \frac{1}{\eta^{s-1}}) \int_{B_{1/\delta}} |\phi_\delta| dx' - (s-1) \int_\varepsilon^\eta \int_{B_{1/\delta}} \frac{|\phi_\delta|}{x_n^s} dx' dx_n}{\int_\varepsilon^\eta \int_{B_{1/\delta}} \frac{|\nabla_{x'}\phi_\delta|}{x_n^\alpha} dx' dx_n + (\frac{1}{\varepsilon^\alpha} + \frac{1}{\eta^\alpha}) \int_{B_{1/\delta}} |\phi_\delta| dx'} \\
&= \frac{K_\delta \int_\varepsilon^\eta x_n^{1-s} dx_n + M_\delta(\varepsilon^{1-s} + \eta^{1-s}) + M_\delta \int_\varepsilon^\eta (x_n^{1-s})' dx_n}{K_\delta \int_\varepsilon^\eta x_n^{-\alpha} dx_n + M_\delta(\varepsilon^{-\alpha} + \eta^{-\alpha})},
\end{aligned}$$

where we have set $K_\delta := \int_{B_{1/\delta}} |\nabla_{x'}\phi_\delta(x')| dx'$ and $M_\delta := \int_{B_{1/\delta}} |\phi_\delta(x')| dx'$. Performing the integration appeared in the last term of the numerator we arrive at

$$\tilde{Q}[u_{\varepsilon,\delta}] = \frac{K_\delta \int_\varepsilon^\eta x_n^{1-s} dx_n + 2M_\delta \eta^{1-s}}{K_\delta \int_\varepsilon^\eta x_n^{-\alpha} dx_n + M_\delta(\varepsilon^{-\alpha} + \eta^{-\alpha})}.$$

By the change of variables $y' = \delta x'$ we obtain $K_\delta = \delta^{2-n} K_1$ and $M_\delta = \delta^{1-n} M_1$. Thus

$$\begin{aligned}
\tilde{Q}[u_{\varepsilon,\delta}] &= \frac{\delta^{2-n} K_1 \int_\varepsilon^\eta x_n^{1-s} dx_n + 2\delta^{1-n} M_1 \eta^{1-s}}{\delta^{2-n} K_1 \int_\varepsilon^\eta x_n^{-\alpha} dx_n + \delta^{1-n} M_1(\varepsilon^{-\alpha} + \eta^{-\alpha})} \\
&= \frac{\delta K_1 \int_\varepsilon^\eta x_n^{1-s} dx_n + 2M_1 \eta^{1-s}}{\frac{\delta K_1}{1-\alpha}(\eta^{1-\alpha} - \varepsilon^{1-\alpha}) + M_1(\varepsilon^{-\alpha} + \eta^{-\alpha})}.
\end{aligned}$$

To proceed we distinguish cases:

- Let $1 < s < 2$. Then

$$\tilde{Q}[u_{\varepsilon,\delta}] = \frac{\frac{\delta K_1}{2-s}(\eta^{2-s} - \varepsilon^{2-s}) + 2M_1\eta^{1-s}}{\frac{\delta K_1}{1-\alpha}(\eta^{1-\alpha} - \varepsilon^{1-\alpha}) + M_1(\varepsilon^{-\alpha} + \eta^{-\alpha})}$$

$$= o_\varepsilon(1).$$

- Now let $s = 2$. Then

$$\tilde{Q}[u_{\varepsilon,\delta}] = \frac{\delta K_1 \log(\eta/\varepsilon) + 2M_1\eta^{-1}}{\frac{\delta K_1}{1-\alpha}(\eta^{1-\alpha} - \varepsilon^{1-\alpha}) + M_1(\varepsilon^{-\alpha} + \eta^{-\alpha})}$$

$$= o_\varepsilon(1).$$

- Finally let $s > 2$. Then

$$\tilde{Q}[u_{\varepsilon,\delta}] = \frac{\frac{\delta K_1}{s-2}(\varepsilon^{2-s} - \eta^{2-s}) + 2M_1\eta^{1-s}}{\frac{\delta K_1}{1-\alpha}(\eta^{1-\alpha} - \varepsilon^{1-\alpha}) + M_1(\varepsilon^{-\alpha} + \eta^{-\alpha})}.$$

We may set $\delta = \varepsilon^{s-2}$ so that $\tilde{Q}[u_{\varepsilon,\delta}] = o_\varepsilon(1)$. ∎

## 2.3 Sets with positive reach

In this subsection we obtain an interpolation inequality between (2.4) and (2.9) via sets with positive reach.

Let $\emptyset \neq K \subsetneq \mathbb{R}^n$ be closed and consider the distance function to $K$ i.e. $d_K : \mathbb{R}^n \to [0, \infty)$ with $d_K(x) = \inf\{|x - y| : y \in K\}$. Denote by $K_1$ the set of points in $\mathbb{R}^n$ which have a unique closest point on $K$, namely $K_1 = \{x \in \mathbb{R}^n : \exists! \, y \in K \text{ such that } d_K(x) = |x - y|\}$.

**Definition 2.15.** The *reach of a point* $x \in K$ is $\operatorname{reach}(K, x) := \sup\{r \geq 0 : B_r(x) \subset K_1\}$. The *reach of the set* $K$ is $\operatorname{reach}(K) := \inf_{x \in K} \operatorname{reach}(K, x)$.

The above definition was introduced in [F] where it was also noted that $K$ is convex if and only if $\operatorname{reach}(K) = \infty$.

**Lemma 2.16.** *Let* $\Omega \subsetneq \mathbb{R}^n$ *be open and set* $h := \operatorname{reach}(\overline{\Omega}) \geq 0$. *Then*

$$(h + d)(-\Delta d) \geq -(n-1) \text{ in } \Omega, \text{ in the sense of distributions,} \tag{2.12}$$

*where* $d \equiv d(x) = \inf\{|x - y| : y \in \mathbb{R}^n \setminus \Omega\}$.

*Proof.* If $h = 0$, this is Lemma 2.1. For $h > 0$ we set $\Omega_h = \{x \in \mathbb{R}^n : d_{\overline{\Omega}}(x) < h\}$. As in the proof of Lemma 2.1, the continuous function $\bar{A} : \mathbb{R}^n \to \mathbb{R}$ defined by $\bar{A}(x) = |x|^2 - d_{\Omega_h^c}^2(x)$ is convex, and thus the distributional Laplacian of $\bar{A}$ is a nonnegative Radon measure on $\mathbb{R}^n$. The result follows since for $x \in \Omega$ we have $d_{\Omega_h^c}(x) = d(x) + h$ (see also [F]-Corollary 4.9), and thus $\Delta \bar{A} = 2(n - 1 - (h + d)\Delta d) \geq 0$ in $\Omega$, in the sense of distributions. ∎

**Theorem 2.17.** *Let $\Omega \subsetneq \mathbb{R}^n$ be open and set $h := \text{reach}(\overline{\Omega})$. Suppose in addition that $R := \sup_{x \in \Omega} d(x) < \infty$. For all $u \in C_c^\infty(\Omega)$ and all $s > \frac{h+nR}{h+R}$, it holds that*

$$\int_\Omega \frac{|\nabla u|}{d^{s-1}} dx \geq \left((s-1)\frac{h}{h+R} + (s-n)\frac{R}{h+R}\right) \int_\Omega \frac{|u|}{d^s} dx. \tag{2.13}$$

*Proof.* Inserting (2.12) to (2.3), we obtain

$$\begin{aligned}
\int_\Omega \frac{|\nabla u|}{d^{s-1}} dx &\geq (s-1)\int_\Omega \frac{|u|}{d^s} dx - (n-1)\int_\Omega \frac{|u|}{d^s}\frac{d}{h+d} dx. \\
&= \int_\Omega \frac{(s-1)h + (s-n)d}{h+d}\frac{|u|}{d^s} dx \\
&\geq \frac{(s-1)h + (s-n)R}{h+R}\int_\Omega \frac{|u|}{d^s} dx,
\end{aligned}$$

where the last inequality follows since $R < \infty$ and $\frac{(s-1)h+(s-n)d}{h+d}$ is decreasing in $d$. ∎

Note that this inequality interpolates between the case of a general open set $\Omega \subsetneq \mathbb{R}^n$, where we have $h = 0$ and the constant becomes $s - n$, and the case of a convex set $\Omega$, where $h = \infty$ and the constant becomes $s - 1$.

# 3  A lower bound on $-\Delta d$ and the role of mean convexity

Before stating our result in this section (Theorem B of the introduction), we gather some additional properties of the distance function to the boundary that will be in use.

From now on $\Omega$ will be a domain, i.e., an open and connected subset of $\mathbb{R}^n$. We will denote by $\Sigma$ the set of points in $\Omega$ which have more than one projection on $\partial\Omega$. If $x \in \Omega \setminus \Sigma$, then $\xi(x)$ will stand for its unique projection on the boundary.

The following Lemma follows from Lemmas 14.16 and 14.17 in [GTr].

**Lemma 3.1.** *Let $\Omega \subset \mathbb{R}^n$ be a domain (possibly unbounded) with boundary of class $\mathcal{C}^2$.*

*(1) If in addition $\Omega$ satisfies a uniform interior sphere condition, then there exists $\delta > 0$ such that $\tilde{\Omega}_\delta := \{x \in \overline{\Omega} : d(x) < \delta\} \subset \overline{\Omega} \setminus \Sigma$ and $d \in C^2(\tilde{\Omega}_\delta)$.*

*(2) $d \in C^2(\overline{\Omega} \setminus \overline{\Sigma})$ and for any $x \in \overline{\Omega} \setminus \overline{\Sigma}$, in terms of a principal coordinate system at $\xi(x) \in \partial\Omega$, it holds that*

$$\begin{aligned}
(i) \quad & \nabla d(x) = -\vec{\nu}(\xi(x)) = (0, ..., 0, 1) \\
(ii) \quad & 1 - \kappa_i(\xi(x))d(x) > 0 \text{ for all } i = 1, ..., n-1 \\
(iii) \quad & [D^2 d(x)] = \text{diag}\left[\frac{-\kappa_1(\xi(x))}{1 - \kappa_1(\xi(x))d(x)}, ..., \frac{-\kappa_{n-1}(\xi(x))}{1 - \kappa_{n-1}(\xi(x))d(x)}, 0\right],
\end{aligned}$$

*where $\vec{\nu}(\xi(x))$ is the unit outer normal at $\xi(x) \in \partial\Omega$, and $\kappa_1(\xi(x)), ..., \kappa_{n-1}(\xi(x))$ are the principal curvatures of $\partial\Omega$ at the point $\xi(x) \in \partial\Omega$.*

**Remark 3.2.** Part (2) of the above Lemma is proved in [GTr] only in $\tilde{\Omega}_\delta$. However, it is also true for the largest open set contained in $\Omega \setminus \Sigma$, i.e. $\Omega \setminus \overline{\Sigma}$ (see for instance [CrM], [LN], [CC], [G]).

Another known, important fact we will need is that domains with boundary of class $C^2$ satisfy $\mathcal{L}^n(\overline{\Sigma}) = 0$. This is proved in [Mnn]-Errata-§5.2 (see also [CrM] where however, only bounded domains are discussed). At last, we shall need the following Lemma for which we add the proof in correspondence to Lemmas 2.1 & 2.16 (see [CS]-Proposition 2.2.2.(ii) & Proposition 1.1.3.(c) and also [Fu]).

**Lemma 3.3.** *Let $\Omega \subsetneq \mathbb{R}^n$ be open. The function $\tilde{A} : \mathbb{R}^n \to \mathbb{R}$ defined by $\tilde{A}(x) = C|x|^2/2 - d(x)$, is convex in any open ball $B \subset\subset \Omega$, for any $C \geq 1/\operatorname{dist}(B, \partial\Omega)$.*

*Proof.* First note that for all $a, b \in \mathbb{R}^n$ with $a \neq 0$, we have

$$|a+b| + |a-b| - 2|a| \leq \frac{|b|^2}{|a|}. \tag{3.1}$$

We choose an open ball $B \subset \Omega$ with $r := \operatorname{dist}(B, \partial\Omega) > 0$, and take $x \in B$. Let $y \in \partial\Omega$ be such that $d(x) = |x - y|$. For any $z \in \mathbb{R}^n$ such that $x + z, x - z \in B$, we get

$$\begin{aligned}
\tilde{A}(x+z) + \tilde{A}(x-z) - 2\tilde{A}(x) &= C|z|^2 - (d(x+z) + d(x-z) - 2d(x)) \\
&\geq C|z|^2 - (|x+z-y| + |x-z-y| - 2|x-y|) \\
\text{(by (3.1) for } a = x - y \text{ and } b = z\text{)} \quad &\geq C|z|^2 - \frac{|z|^2}{|x-y|} \\
&\geq (C - 1/r)|z|^2.
\end{aligned}$$

Since $\tilde{A}(x)$ is also continuous, we obtain that $\tilde{A}(x)$ is convex in $B$ for any $C \geq 1/r$. ∎

To state our main result in this subsection, we denote by $\mathcal{H}(y) := \frac{1}{n-1}\sum_{i=1}^{n-1}\kappa_i(y)$ the mean curvature of $\partial\Omega$ at the point $y \in \partial\Omega$.

**Theorem 3.4.** *Let $\Omega \subset \mathbb{R}^n$ be a domain with boundary of class $C^2$ satisfying a uniform interior sphere condition. Then $\mu := (-\Delta d)dx$ is a signed Radon measure on $\Omega$. Let $\mu = \mu_{ac} + \mu_s$ be the Lebesgue decomposition of $\mu$ with respect to $\mathcal{L}^n$, i.e., $\mu_{ac} \ll \mathcal{L}^n$ and $\mu_s \perp \mathcal{L}^n$. Then $\mu_s \geq 0$ in $\Omega$, and $\mu_{ac} \geq (n-1)\underline{\mathcal{H}}dx$ a.e. in $\Omega$, where $\underline{\mathcal{H}} := \inf_{y \in \partial\Omega}\mathcal{H}(y)$.*

*Proof.* Letting $\delta$ be as in Lemma 3.1(1), we set $\Omega_\delta = \{x \in \Omega : d(x) < \delta\}$. Then $-\Delta d$ is a continuous function on $\Omega_\delta$ and so $\mu^0 := (-\Delta d)dx$ is a signed Radon measure on $\Omega_\delta$, absolutely continuous with respect to $\mathcal{L}^n$.

Next, let $\{B_i\}_{i\geq 1}$ be a cover of the set $\Omega \setminus \Omega_\delta$, comprised of open balls $B_i$ for which $\operatorname{dist}(B_i, \partial\Omega) > \delta/2$ for all $i \geq 1$. According to Lemma 3.3, the function $\tilde{A}(x) := |x|^2/\delta - d(x)$ is convex in each $B_i$. From [EvG]-§6.3-Theorem 2, we deduce that there exist nonnegative Radon measures $\{\nu^i\}_{i\geq 1}$, respectively on $\{B_i\}_{i\geq 1}$, such that

$$\int_{B_i} \phi \Delta \tilde{A} dx = \int_{B_i} \phi d\nu^i,$$

for all $\phi \in C_c^\infty(B_i)$. Since $\Delta \tilde{A} = 2n/\delta - \Delta d$ in the sense of distributions, we get

$$\int_{B_i} \phi(-\Delta d)\mathrm{d}x = \int_{B_i} \phi \mathrm{d}\nu^i - \frac{2n}{\delta}\int_{B_i}\phi \mathrm{d}x, \tag{3.2}$$

for all $\phi \in C_c^\infty(B_i)$, and thus $\mu^i := (-\Delta d)\mathrm{d}x = \nu^i - \frac{2n}{\delta}\mathrm{d}x$ is a signed Radon measure on $B_i$.

Let $\{\eta_i\}_{i\geq 1}$ be a $C^\infty$ partition of unity subordinated to the open covering $\{B_i\}_{i\geq 1}$ of $\Omega\setminus\Omega_\delta$, i.e.

$$\eta_i \in C_c^\infty(B_i), \quad 0 \leq \eta_i(x) \leq 1 \text{ in } B_i \quad \text{and} \quad \sum_{i=1}^\infty \eta_i(x) = 1 \text{ in } \Omega\setminus\Omega_\delta.$$

Further, for $x \in \Omega$ define $\eta_0(x) = 1 - \sum_{i=1}^\infty \eta_i(x)$. We then have

$$\mathrm{sprt}\,\eta_0 \subset \Omega_\delta, \quad \eta_0(x) = 1 \text{ in } \Omega_{\delta/2} \quad \text{and} \quad \sum_{i=0}^\infty \eta_i(x) = 1 \text{ in } \Omega.$$

We will now show that $\mu := \sum_{i=0}^\infty \eta_i \mu^i$ is a well defined signed Radon measure on $\Omega$, and $\mu = (-\Delta d)\mathrm{d}x$. To this end, for any $\phi \in C_c^\infty(\Omega)$ we have

$$\begin{aligned}
\int_\Omega \phi(-\Delta d)\mathrm{d}x &= \sum_{i=0}^\infty \int_\Omega \phi\eta_i(-\Delta d)\mathrm{d}x \\
(\text{by (3.2)}) &= \int_\Omega \phi\eta_0 \mathrm{d}\mu^0 + \sum_{i=1}^\infty \left(\int_\Omega \phi\eta_i \mathrm{d}\nu^i - \frac{2n}{\delta}\int_\Omega \phi\eta_i \mathrm{d}x\right) \\
&= \int_\Omega \phi\eta_0 \mathrm{d}\mu^0 + \int_\Omega \phi\sum_{i=1}^\infty \eta_i \mathrm{d}\nu^i - \frac{2n}{\delta}\int_\Omega \phi\sum_{i=1}^\infty \eta_i \mathrm{d}x \\
&= \int_\Omega \phi\eta_0 \mathrm{d}\mu^0 + \int_\Omega \phi\sum_{i=1}^\infty \eta_i \mathrm{d}\mu^i \\
&= \int_\Omega \phi \mathrm{d}\mu,
\end{aligned}$$

where the middle equality follows since $\nu^i$ are positive Radon measures and thus $\sum_{i=0}^m \eta_i \nu^i$ is increasing in $m$ (see [EvG]-Section 1.9).

Next, by the Lebesgue Decomposition Theorem ([EvG]-§1.3-Theorem 3), $\mu = \mu_{ac} + \mu_s$ where

$$\mu_s = \sum_{i=0}^\infty \eta_i \mu_s^i = \sum_{i=1}^\infty \eta_i \mu_s^i = \sum_{i=1}^\infty \eta_i \nu_s^i \geq 0,$$

since $\mu^i = \nu^i - \frac{2n}{\delta}\mathrm{d}x$ and $\nu^i$ are nonnegative. Finally, from Lemma 3.1-(2) we get

$$\begin{aligned}
-\Delta d(x) &= \sum_{i=1}^{n-1} \frac{\kappa_i(\xi(x))}{1 - \kappa_i(\xi(x))d(x)} \\
&\geq \sum_{i=1}^{n-1} \kappa_i(\xi(x)) \\
&= (n-1)\mathcal{H}(\xi(x)) \\
&\geq (n-1)\underline{\mathcal{H}}, \quad \forall x \in \Omega \setminus \overline{\Sigma}.
\end{aligned}$$

Now by Lemma 3.1-(2), $-\Delta d$ is a continuous function on $\Omega \setminus \overline{\Sigma}$ and so

$$\mu_{ac} = (-\Delta d)\mathrm{d}x \geq (n-1)\underline{\mathcal{H}}\mathrm{d}x \text{ in } \Omega \setminus \overline{\Sigma}.$$

Recalling that $\mathcal{L}^n(\overline{\Sigma}) = 0$ when $\partial\Omega \in \mathcal{C}^2$ and since $\Omega = (\Omega \setminus \overline{\Sigma}) \cup \overline{\Sigma}$, we conclude $\mu_{ac} \geq (n-1)\underline{\mathcal{H}}\mathrm{d}x$ a.e. in $\Omega$. ∎

**Definition 3.5.** A domain $\Omega$ with boundary of class $\mathcal{C}^2$ is said to be *mean convex* if $\mathcal{H}(y) \geq 0$ for all $y \in \partial\Omega$.

Theorem 3.4 along with Lemma 3.1 provides us a characterization of mean convexity in terms of the distance function for sufficiently smooth domains. More precisely, we have the following

**Corollary 3.6.** *Let $\Omega$ be a domain with boundary of class $\mathcal{C}^2$ satisfying a uniform interior sphere condition. Then $\Omega$ is mean convex if and only if condition (C) holds, i.e., $-\Delta d \geq 0$ holds in $\Omega$, in the sense of distributions.*

**Remark 3.7.** The resulting lower bound $-(\Delta d)\mathrm{d}x \geq (n-1)\underline{\mathcal{H}}\mathrm{d}x$, is optimal. To see this, assume first that $\Omega$ is bounded and choose a point $y_0 \in \partial\Omega$ such that $\mathcal{H}(y_0) = \underline{\mathcal{H}}$. Pick $0 \leq \phi_\delta \in C_c^\infty(\Omega)$, such that $\mathrm{sprt}\{\phi_\delta\} \subset B_\delta(y_0) \cap \Omega_\delta$, where $\delta > 0$, small. For sufficiently small $\delta$, as $x \in \Omega_\delta$ approaches $y_0$ we have $-\Delta d(x) = (n-1)\underline{\mathcal{H}} + O(d(x))$. Thus, as $\delta \downarrow 0$ we have $-\Delta d(x) = (n-1)\underline{\mathcal{H}} + o_\delta(1)$ for all $x \in B_\delta(y_0) \cap \Omega_\delta$, and so

$$\begin{aligned}
\inf_{0 \leq \phi \in C_c^\infty(\Omega) \setminus \{0\}} \frac{\int_\Omega \phi(-\Delta d)\mathrm{d}x}{\int_\Omega \phi\mathrm{d}x} &\leq \frac{\int_{B_\delta(y_0) \cap \Omega_\delta} \phi_\delta(-\Delta d)\mathrm{d}x}{\int_{B_\delta(y_0) \cap \Omega_\delta} \phi_\delta \mathrm{d}x} \\
&= (n-1)\underline{\mathcal{H}} + o_\delta(1).
\end{aligned}$$

If $\Omega$ is unbounded, we may consider a sequence $\{y_n\} \subset \partial\Omega$ converging to $y_0$, and repeat the above argument for any such point, to obtain

$$\inf_{0 \leq \phi \in C_c^\infty(\Omega) \setminus \{0\}} \frac{\int_\Omega \phi(-\Delta d)\mathrm{d}x}{\int_\Omega \phi\mathrm{d}x} \leq (n-1)\mathcal{H}(y_n) + o_\delta(1).$$

Since $\mathcal{H}(y)$ is a continuous function on $\partial\Omega$, we end up by letting $n \to \infty$.

# 4 Proof of Theorem A and Theorem D

Let $\Omega$ be a domain satisfying property (C). We define the quotient

$$Q_\beta[u] := \frac{\int_\Omega \frac{|\nabla u|}{d^{s-1}} dx - (s-1) \int_\Omega \frac{|u|}{d^s} dx}{\int_\Omega \frac{|u|}{d^{s-\beta}} dx}; \quad s > 1, \tag{4.1}$$

and we consider the following minimization problem

$$\mathcal{B}_\beta(\Omega) := \inf\{Q_\beta[u] : u \in C_c^\infty(\Omega) \setminus \{0\}\}; \quad 0 < \beta \leq s - 1.$$

The next Proposition shows that the essential range for $\beta$ is smaller.

**Proposition 4.1.** *Let $\Omega$ be a domain with boundary of class $\mathcal{C}^2$ satisfying property $(C)$. If $s \geq 2$ then $\mathcal{B}_\beta(\Omega) = 0$ for all $0 < \beta < 1$. If $1 < s < 2$ then $\mathcal{B}_\beta(\Omega) = 0$ for all $0 < \beta \leq s - 1$.*

*Proof.* For small $\delta > 0$, let $\Omega_\delta := \{x \in \Omega : d(x) < \delta\}$ and $\Omega_\delta^c = \Omega \setminus \Omega_\delta$. We test $(4.1)$ with $u_\delta(x) = \chi_{\Omega_\delta^c}(x)\phi(x)$, where $\phi \in C_c^\infty(B_\varepsilon(y_0))$ for a fixed $y_0 \in \partial\Omega$ and sufficiently small $\varepsilon$, satisfying $\varepsilon > 3\delta$. We may suppose in addition that $0 \leq \phi \leq 1$ in $B_\varepsilon(y_0)$, $\phi \equiv 1$ in $B_{\varepsilon/2}(y_0)$ and $|\nabla \phi| \leq 1/\varepsilon$. This function is not in $C_c^\infty(\Omega)$, but since it is in $BV(\Omega)$ we can mollify the characteristic function so that the calculations below to hold in the limit. The distributional gradient of $u_\delta$ is $\nabla u_\delta = \chi_{\Omega_\delta^c} \nabla \phi - \vec{\nu} \phi \delta_{\partial\Omega_\delta^c}$, where $\vec{\nu}$ is the outward pointing, unit normal vector field along $\partial\Omega_\delta^c$, and $\delta_{\partial\Omega_\delta^c}$ is the Dirac measure on $\partial\Omega_\delta^c$. Moreover, the total variation of $\nabla u_\delta$ is $|\nabla u_\delta| = \chi_{\Omega_\delta^c}|\nabla \phi| + \phi \delta_{\partial\Omega_\delta^c}$. Since $\partial\Omega_\delta^c = \{x \in \Omega : d(x) = \delta\}$, we obtain

$$Q_\beta[u_\delta] = \frac{\int_{\Omega_\delta^c} |\nabla \phi| d^{1-s} dx + \delta^{1-s} \int_{\partial\Omega_\delta^c} \phi dS_x - (s-1) \int_{\Omega_\delta^c} \phi d^{-s} dx}{\int_{\Omega_\delta^c} \phi d^{\beta-s} dx}. \tag{4.2}$$

Using the fact that $|\nabla d(x)| = 1$ for a.e. $x \in \Omega$, we may perform an integration by parts in the last term of the numerator as follows

$$\begin{aligned}
(s-1) \int_{\Omega_\delta^c} \phi d^{-s} dx &= -\int_{\Omega_\delta^c} \phi \nabla d \cdot \nabla d^{1-s} dx \\
&= \int_{\Omega_\delta^c} [\nabla \phi \cdot \nabla d] d^{1-s} dx + \int_{\Omega_\delta^c} \phi d^{1-s} \Delta d\, dx - \delta^{1-s} \int_{\partial\Omega_\delta^c} \phi \nabla d \cdot \vec{\nu} dS_x.
\end{aligned}$$

Since $\nabla d$ is the inner unit normal to $\partial\Omega$, we have $\nabla d \cdot \vec{\nu} = -1$ and substituting the above equality in $(4.2)$, the surface integrals will be canceled to get

$$Q_\beta[u_\delta] = \frac{\int_{\Omega_\delta^c} [|\nabla \phi| - \nabla \phi \cdot \nabla d] d^{1-s} dx + \int_{\Omega_\delta^c} \phi d^{1-s}(-\Delta d) dx}{\int_{\Omega_\delta^c} \phi d^{\beta-s} dx}.$$

By the fact that $-\Delta d(x) \leq c$ for all $x \in \Omega_\delta^c \cap B_\varepsilon$, and by the properties we imposed on $\phi$, we get

$$Q_\beta[u_\delta] \leq \frac{\frac{2}{\varepsilon}\int_{\Omega_\delta^c \cap B_\varepsilon} d^{1-s}\mathrm{d}x + c\int_{\Omega_\delta^c \cap B_\varepsilon} d^{1-s}\mathrm{d}x}{\int_{\Omega_\delta^c \cap B_{\varepsilon/2}} d^{\beta-s}\mathrm{d}x}$$

$$= c(\varepsilon)\frac{\int_{\Omega_\delta^c \cap B_\varepsilon} d^{1-s}\mathrm{d}x}{\int_{\Omega_\delta^c \cap B_{\varepsilon/2}} d^{\beta-s}\mathrm{d}x}$$

$$=: c(\varepsilon)\frac{N(\delta)}{D(\delta)}.$$

Using now the co-area formula we compute

$$N(\delta) = \int_\delta^\varepsilon r^{1-s}\int_{\{x\in\Omega_\delta^c \cap B_\varepsilon: d(x)=r\}} \mathrm{d}S_x \mathrm{d}r$$

$$\leq c_1(\varepsilon)\int_\delta^\varepsilon r^{1-s}\mathrm{d}r,$$

where $c_1(\varepsilon) = \max_{r\in[0,\varepsilon]}|\{x \in \Omega_\delta^c \cap B_\varepsilon : d(x) = r\}|$. Also,

$$D(\delta) = \int_\delta^{\varepsilon/2} r^{\beta-s}\int_{\{x\in\Omega_\delta^c \cap B_{\varepsilon/2}: d(x)=r\}} \mathrm{d}S_x \mathrm{d}r$$

$$\geq \int_\delta^{\varepsilon/3} r^{\beta-s}\int_{\{x\in\Omega_\delta^c \cap B_{\varepsilon/2}: d(x)=r\}} \mathrm{d}S_x \mathrm{d}r$$

$$\geq c_2(\varepsilon)\int_\delta^{\varepsilon/3} r^{\beta-s}\mathrm{d}r,$$

where $c_2(\varepsilon) = \min_{r\in[0,\varepsilon/3]}|\{x \in \Omega_\delta^c \cap B_{\varepsilon/2} : d(x) = r\}|$. A direct computation reveals that if $s \geq 2$ then $Q_\beta[u_\delta] \leq o_\delta(1)$ for all $0 < \beta < 1$, and also if $1 < s < 2$ then $Q_\beta[u_\delta] \leq o_\delta(1)$ for all $0 < \beta \leq s - 1$. ∎

## 4.1 Lower and upper estimates for $\mathcal{B}_1(\Omega)$

In this subsection we obtain upper and lower estimates for $\mathcal{B}_1(\Omega)$. In particular we prove Theorem A and the optimality in Theorem D of the introduction.

**Theorem 4.2** (**Lower estimate**). *Let $\Omega$ be a domain with boundary of class $C^2$ satisfying a uniform interior sphere condition. If $s \geq 1$ then*

$$\mathcal{B}_1(\Omega) \geq (n-1)\underline{\mathcal{H}}, \tag{4.3}$$

*where $\underline{\mathcal{H}}$ is the infimum of the mean curvature of $\partial\Omega$.*

*Proof.* The estimate follows directly from (2.3) using Theorem 3.4. ∎

**Remark 4.3.** By Theorem 2.9, if condition (C) is satisfied, then the first term in (2.3) is sharp. The passage from (2.3) to inequality (4.3) via Theorem 3.4, is also sharp, i.e. the constant $(n-1)\underline{\mathcal{H}}$ in the inequality

$$\int_\Omega \frac{|u|}{d^{s-1}}(-\Delta d)\mathrm{d}x \geq (n-1)\underline{\mathcal{H}} \int_\Omega \frac{|u|}{d^{s-1}}\mathrm{d}x, \quad \forall u \in C_c^\infty(\Omega),$$

is optimal. To see this, set $v = d^{1-s}|u|$, to get

$$\inf_{u \in C_c^\infty(\Omega)\setminus\{0\}} \frac{\int_\Omega \frac{|u|}{d^{s-1}}(-\Delta d)\mathrm{d}x}{\int_\Omega \frac{|u|}{d^{s-1}}\mathrm{d}x} \leq \inf_{0 \leq v \in C_c^\infty(\Omega)\setminus\{0\}} \frac{\int_\Omega v(-\Delta d)\mathrm{d}x}{\int_\Omega v\mathrm{d}x}$$
$$\leq (n-1)\underline{\mathcal{H}} + o_\delta(1),$$

by Remark 3.7.

We next present upper bounds. We begin with an upper bound which, although not sharp enough for our problem, it is of independent interest.

**Definition 4.4.** The *Cheeger constant* $h(\Omega)$ of a bounded domain $\Omega$ with piecewise $\mathcal{C}^1$ boundary, is defined by $h(\Omega) := \inf_\omega \frac{|\partial \omega|}{|\omega|}$, where the infimum is taken over all sub-domains $\omega \subset\subset \Omega$ with piecewise $\mathcal{C}^1$ boundary.

For existence of minimizers, uniqueness and regularity results concerning the Cheeger constant, we refer to [FrK] and references therein (especially [StrZ]).

**Proposition 4.5.** *Let $\Omega$ be a bounded domain with piecewise $\mathcal{C}^1$ boundary such that condition (C) holds. For all $s \geq 1$, we have $\mathcal{B}_1(\Omega) \leq h(\Omega)$.*

*Proof.* Take $\omega \subset\subset \Omega$ with piecewise $\mathcal{C}^1$ boundary and let $u_\omega(x) = (d(x))^{s-1}\chi_\omega(x)$. The distributional gradient and the total variation of this $BV(\Omega)$ function, are respectively, $\nabla u_\omega = (s-1)d^{s-2}\chi_\omega \nabla d - \vec{\nu}d^{s-1}\delta_{\partial \omega}$ and $|\nabla u_\omega| = (s-1)d^{s-2}\chi_\omega + d^{s-1}\delta_{\partial \omega}$, where $\vec{\nu}$ is the outward pointing, unit normal vector field along $\partial \omega$, and $\delta_{\partial \omega}$ is the uniform Dirac measure on $\partial \omega$. We test (4.1) with $u_\omega$ to get

$$Q_1[u_\omega] = \frac{(s-1)\int_\omega d^{-1}\mathrm{d}x + \int_{\partial\omega} \mathrm{d}S_x - (s-1)\int_\omega d^{-1}\mathrm{d}x}{\int_\omega \mathrm{d}x} = \frac{|\partial \omega|}{|\omega|}.$$

In particular $h(\Omega) = \inf_\omega Q_1[u_\omega]$. By the standard $C_c^\infty$ approximation of the characteristic function of the domain $\omega$, we obtain $\mathcal{B}_1(\Omega) \leq \frac{|\partial \omega|}{|\omega|}$ and thus $\mathcal{B}_1(\Omega) \leq h(\Omega)$. ∎

From Theorem 4.2 and Proposition 4.5 for $s = 1$, we conclude the following

**Corollary 4.6.** *If $\Omega$ is a strictly mean convex, bounded domain with boundary of class $\mathcal{C}^2$, it holds that $h(\Omega) \geq (n-1)\underline{\mathcal{H}}$.*

**Remark.** In [AltC] it is proved that a bounded convex domain $\Omega$ is a self-minimizer of $h(\Omega)$, *if and only if* it belongs to the class $\mathcal{C}^{1,1}$ and also the stronger estimate $h(\Omega) \geq (n-1)\overline{\mathcal{H}}$ holds. Here $\overline{\mathcal{H}}$ is the essential supremum of the mean curvature of the boundary (the last being defined in the almost everywhere sense since $\partial\Omega \in \mathcal{C}^{1,1}$).

**Remark.** By Corollary 4.6, if $\Omega$ is a bounded strictly mean convex domain with boundary of class $\mathcal{C}^2$, then

$$\frac{|\partial\Omega|}{|\Omega|} \geq (n-1)\underline{\mathcal{H}}. \tag{4.4}$$

For bounded convex domains with boundary of class $\mathcal{C}^2$, this follows by one of the Minkowski quadratic inequalities for cross-sectional measures (see [BZ] - eq(16), pg 144). It states that

$$\frac{|\partial\Omega|}{|\Omega|} \geq \frac{n}{|\partial\Omega|} \int_{\partial\Omega} \mathcal{H}(y)\mathrm{d}S.$$

Thus we have $|\partial\Omega|/|\Omega| \geq n\underline{\mathcal{H}}$ from which (4.4) follows. This remark is taken from [GN], where one can also find an application of (4.4).

The following result states a more useful upper bound for $\mathcal{B}_1(\Omega)$. It will be combined with Theorem 4.2 to give the best possible constant for special geometries.

**Theorem 4.7.** *Let $\Omega$ be a domain with boundary of class $\mathcal{C}^2$ satisfying a uniform interior sphere condition. If $s \geq 2$ then for all $\phi \in C_c^1(\partial\Omega)$,*

$$\mathcal{B}_1(\Omega) \leq (n-1)\frac{\int_{\partial\Omega}|\phi(y)|\mathcal{H}(y)\mathrm{d}S}{\int_{\partial\Omega}|\phi(y)|\mathrm{d}S} + \frac{\int_{\partial\Omega}|\nabla\phi(y)|\mathrm{d}S}{\int_{\partial\Omega}|\phi(y)|\mathrm{d}S},$$

*where $\mathcal{H}(y)$ is the mean curvature at the point $y \in \partial\Omega$.*

*Proof.* Let $\delta > 0$ such that for all $x \in \tilde{\Omega}_\delta := \{x \in \overline{\Omega} : d(x) < \delta\}$ there exists a unique point

$$\xi \equiv \xi(x) = x - d(x)\nabla d(x) \in \partial\Omega \tag{4.5}$$

with $d(x) = |x - \xi|$. For any $t \in [0, \delta]$ the surface area element of $\partial\Omega_t^c = \{x \in \Omega : d(x) = t\}$ is given by

$$\mathrm{d}S_t = (1 - \kappa_1 t)...(1 - \kappa_{n-1} t)\mathrm{d}S = (1 - (n-1)t\mathcal{H} + O(t^2))\mathrm{d}S, \tag{4.6}$$

where $\kappa_1, ..., \kappa_{n-1}$, are the principal curvatures of $\partial\Omega$, $\mathrm{d}S$ is the surface area element of $\partial\Omega$ and $\mathcal{H}$ is the mean curvature of $\partial\Omega$, (see [S]-§13.5 & 13.6). Now let $0 < \varepsilon < \delta$ and chose $\phi \in C_c^1(\partial\Omega)$. We test (4.1) with $u_\varepsilon(x) = \chi_{\Omega_\varepsilon^c \setminus \Omega_\delta^c}(x)\phi(\xi(x))$, $\xi(x)$ as in (4.5), and then we will check the limit as $\varepsilon \downarrow 0$. The distributional gradient of $u_\varepsilon$, is $\nabla u_\varepsilon = (\vec{\nu}_\delta \delta_{\partial\Omega_\delta^c} - \vec{\nu}_\varepsilon \delta_{\partial\Omega_\varepsilon^c})\phi(\xi) +$ $\chi_{\Omega_\varepsilon^c\setminus\Omega_\delta^c}\nabla_x\phi(\xi)$, where $\vec{\nu}_\delta, \vec{\nu}_\varepsilon$ are respectively, the outward pointing unit normal vector fields along $\partial\Omega_\delta^c, \partial\Omega_\varepsilon^c$. Its total variation is $|\nabla u_\varepsilon| = (\delta_{\partial\Omega_\delta^c} + \delta_{\partial\Omega_\varepsilon^c})|\phi(\xi)| + \chi_{\Omega_\varepsilon^c\setminus\Omega_\delta^c}|\nabla_x\phi(\xi)|$. Thus

$$\int_\Omega \frac{|\nabla u_\varepsilon|}{d^{s-1}}\mathrm{d}x = \delta^{1-s}\int_{\partial\Omega_\delta^c}|\phi(\xi)|\mathrm{d}S_\delta + \varepsilon^{1-s}\int_{\partial\Omega_\varepsilon^c}|\phi(\xi)|\mathrm{d}S_\varepsilon + \int_{\Omega_\varepsilon^c\setminus\Omega_\delta^c}\frac{|\nabla_x\phi(\xi)|}{d^{s-1}}\mathrm{d}x. \quad (4.7)$$

The first integral on the right-hand side of (4.7) is a constant since we will keep $\delta$ fixed. We perform the change of variables $y = \xi(x)$ in the second integral. Using (4.6) we have

$$\begin{aligned}
\varepsilon^{1-s}\int_{\partial\Omega_\varepsilon^c}|\phi(\xi)|\mathrm{d}S_\varepsilon &= \varepsilon^{1-s}\int_{\partial\Omega}|\phi(y)|(1 - (n-1)\varepsilon\mathcal{H}(y) + O(\varepsilon^2))\mathrm{d}S \\
&= \varepsilon^{1-s}M - (n-1)\varepsilon^{2-s}M_\mathcal{H} + O(\varepsilon^{3-s}),
\end{aligned} \quad (4.8)$$

where $M := \int_{\partial\Omega}|\phi|\mathrm{d}S$ and $M_\mathcal{H} := \int_{\partial\Omega}|\phi|\mathcal{H}\mathrm{d}S$. Using the co-area formula, the third term on the right-hand side of (4.7) is written as follows

$$\int_{\Omega_\varepsilon^c\setminus\Omega_\delta^c}\frac{|\nabla_x\phi(\xi)|}{d^{s-1}}\mathrm{d}x = \int_\varepsilon^\delta t^{1-s}\int_{\partial\Omega_t^c}|\nabla_x\phi(\xi)|\mathrm{d}S_t\mathrm{d}t. \quad (4.9)$$

From (4.5) we have $\xi_i(x) = x_i - d(x)\frac{\partial}{\partial x_i}(d(x))$ and thus by Lemma 3.1-(c) we compute

$$\begin{aligned}
\nabla_x\phi(\xi) &= \Big(\sum_{i=1}^n \phi_{\xi_i}(\xi)\frac{\partial\xi_i}{\partial x_1}, ..., \sum_{i=1}^n \phi_{\xi_i}(\xi)\frac{\partial\xi_i}{\partial x_n}\Big) \\
&= \Big(\frac{\phi_{\xi_1}(\xi)}{1-\kappa_1 d}, ..., \frac{\phi_{\xi_{n-1}}(\xi)}{1-\kappa_{n-1}d}, 0\Big).
\end{aligned}$$

Thus, (4.9) becomes

$$\begin{aligned}
\int_{\Omega_\varepsilon^c\setminus\Omega_\delta^c}\frac{|\nabla_x\phi(\xi)|}{d^{s-1}}\mathrm{d}x &= \int_\varepsilon^\delta t^{1-s}\int_{\partial\Omega_t^c}\Big(\sum_{i=1}^{n-1}\Big(\frac{\phi_{y_i}}{1-\kappa_i t}\Big)^2\Big)^{1/2}\mathrm{d}S_t\mathrm{d}t \\
&= \int_\varepsilon^\delta t^{1-s}\int_{\partial\Omega}\Big(\sum_{i=1}^{n-1}\Big(\phi_{y_i}\prod_{j=1,j\neq i}^{n-1}(1-\kappa_j t)\Big)^2\Big)^{1/2}\mathrm{d}S\mathrm{d}t,
\end{aligned}$$

where we have changed variables by $y = \xi(x)$ in the last inequality. Expanding the product as in (4.6), we get

$$\begin{aligned}
\int_{\Omega_\varepsilon^c\setminus\Omega_\delta^c}\frac{|\nabla_x\phi(\xi)|}{d^{s-1}}\mathrm{d}x &\leq \int_\varepsilon^\delta t^{1-s}\int_{\partial\Omega}\Big(\sum_{i=1}^{n-1}\phi_{y_i}^2\Big(1 - [(n-1)\mathcal{H}-\kappa_i]t + c_1 t^2\Big)^2\Big)^{1/2}\mathrm{d}S\mathrm{d}t \\
&\leq K\int_\varepsilon^\delta t^{1-s}\mathrm{d}t + c_2\int_\varepsilon^\delta t^{2-s}\mathrm{d}t,
\end{aligned} \quad (4.10)$$

for some $c_1, c_2 \geq 0$, where $K := \int_{\partial\Omega} |\nabla\phi| \mathrm{d}S$. Next, using co-area formula and the same change of variables, we get

$$(s-1)\int_\Omega \frac{|u_\varepsilon|}{d^s}\mathrm{d}x = (s-1)\int_\varepsilon^\delta t^{-s}\int_{\partial\Omega_t^c} |\phi(\xi)|\mathrm{d}S_t\mathrm{d}t$$

$$\geq (s-1)\int_\varepsilon^\delta t^{-s}\int_{\partial\Omega} |\phi(y)|[1-(n-1)t\mathcal{H}(y)+c_3 t^2]\mathrm{d}S\mathrm{d}t$$

$$= M\varepsilon^{1-s} - (s-1)(n-1)M_\mathcal{H}\int_\varepsilon^\delta t^{1-s}\mathrm{d}t + c_4\int_\varepsilon^\delta t^{2-s}\mathrm{d}t, \quad (4.11)$$

for some $c_3, c_4 \in \mathbb{R}$, and similarly

$$\int_\Omega \frac{|u_\varepsilon|}{d^{s-\beta}}\mathrm{d}x \geq M\int_\varepsilon^\delta t^{\beta-s}\mathrm{d}t - (n-1)M_\mathcal{H}\int_\varepsilon^\delta t^{1+\beta-s}\mathrm{d}t + c_5\int_\varepsilon^\delta t^{2+\beta-s}\mathrm{d}t, \quad (4.12)$$

for some $c_5 \in \mathbb{R}$. Thus inserting (4.8), (4.10), (4.11) into (4.7), and by (4.12) for $\beta = 1$, we get

$$Q_\beta[u_\varepsilon] \leq \frac{(n-1)M_\mathcal{H}[(s-1)\int_\varepsilon^\delta t^{1-s}\mathrm{d}t - \varepsilon^{2-s}] + K\int_\varepsilon^\delta t^{1-s}\mathrm{d}t + c_6\int_\varepsilon^\delta t^{2-s}\mathrm{d}t}{M\int_\varepsilon^\delta t^{\beta-s}\mathrm{d}t - (n-1)M_\mathcal{H}\int_\varepsilon^\delta t^{1+\beta-s}\mathrm{d}t + c_5\int_\varepsilon^\delta t^{2+\beta-s}\mathrm{d}t}, \quad (4.13)$$

for some $c_6 \in \mathbb{R}$. If $s = 2$, then

$$Q_1[u_\varepsilon] \leq \frac{((n-1)M_\mathcal{H} + K)\log(\delta/\varepsilon) + O_\varepsilon(1)}{M\log(\delta/\varepsilon) + O_\varepsilon(1)},$$

while if $s > 2$, then

$$Q_1[u_\varepsilon] \leq \frac{\frac{1}{s-2}((n-1)M_\mathcal{H} + K)\varepsilon^{2-s} + c_7\int_\varepsilon^\delta t^{2-s}\mathrm{d}t}{\frac{1}{s-2}M\varepsilon^{2-s} - (n-1)M_\mathcal{H}\int_\varepsilon^\delta t^{2-s}\mathrm{d}t + c_8\int_\varepsilon^\delta t^{3-s}\mathrm{d}t},$$

for some $c_7, c_8 \in \mathbb{R}$. In any case, letting $\varepsilon \downarrow 0$ we deduce $\mathcal{B}_1(\Omega) \leq \frac{(n-1)M_\mathcal{H}+K}{M}$. ∎

An immediate consequence is

**Corollary 4.8** (**Upper estimate**). *Let $\Omega$ be a bounded domain with boundary of class $\mathcal{C}^2$. If $s \geq 2$ then*

$$\mathcal{B}_1(\Omega) \leq \frac{n-1}{|\partial\Omega|}\int_{\partial\Omega} \mathcal{H}(y)\mathrm{d}S$$

*where $\mathcal{H}(y)$ is the mean curvature at the point $y \in \partial\Omega$.*

*Proof.* Since $\Omega$ is bounded we can chose $\varphi \equiv 1$ in the above Theorem. ∎

The proof of Theorem A follows from Proposition 4.1, Theorem 4.2 and Corollary 4.8.

**Example 4.9** (**Ball**). Let $B_R$ be a ball of radius $R$. By Theorem 4.2 we have $\mathcal{B}_1(B_R) \geq \frac{n-1}{R}$, and by Corollary 4.8, $\mathcal{B}_1(B_R) \leq \frac{n-1}{R}$. We conclude that if $s \geq 2$, then $\mathcal{B}_1(B_R) = \frac{n-1}{R}$. See §5.

**Example 4.10** (**Infinite strip: proof of the optimality in Theorem D**). Let $S_R = \{x = (x', x_n) : x' \in \mathbb{R}^{n-1}, 0 < x_n < 2R\}$. If $s \geq 2$, then combining Theorem 4.2 and Theorem 4.7 we can prove that $\mathcal{B}_1(S_R) = 0$. In fact we have $\mathcal{B}_\beta(S_R) = 0$ for any $1 < \beta \leq s - 1$ and in particular we will prove that if $\gamma = 1$, there is not positive constant $C$ such that (2.10) holds for $\gamma = 1$. To see this, pick any $\phi \equiv \phi(x') \in C_c^1(\mathbb{R}^{n-1})$ such that $\mathrm{sprt}\{\phi\} \subset B_1 \subset \mathbb{R}^{n-1}$, where $B_1$ is the open ball in $\mathbb{R}^{n-1}$ with radius 1, centered at $0'$. Let $\eta > 0$ and set $\phi_\eta \equiv \phi_\eta(x') := \phi(\eta x')$. Note that $\mathrm{sprt}\{\phi_\eta\} \subset B_{1/\eta}$. Let also $0 < \varepsilon < \delta$ for some fixed $\delta \leq R$ (so that $d(x) = x_n$). The quotient corresponding to (2.10) is

$$Q_\gamma[u] = \frac{\int_{S_R} \frac{|\nabla u|}{d^{s-1}} dx - (s-1) \int_{S_R} \frac{|u|}{d^s} dx}{\int_{S_R} \frac{|u|}{d} X^\gamma(\frac{d}{R}) dx} \tag{4.14}$$

As in the proof of Theorem 4.7 we test (4.14) with $u_{\varepsilon,\eta}(x) := \chi_{(\varepsilon,\delta)}(x_n) \phi_\eta(x')$, to arrive at

$$Q_\gamma[u_{\varepsilon,\eta}] = \frac{K_\eta \int_\varepsilon^\delta x_n^{1-s} dx_n + 2M_\eta \delta^{1-s}}{M_\eta \int_\varepsilon^\delta x_n^{-1} X^\gamma(x_n/R) dx_n},$$

where we have set $M_\eta := \int_{B_{1/\eta}} |\phi_\eta(x')| dx'$ and $K_\eta := \int_{B_{1/\eta}} |\nabla_{x'} \phi_\eta(x')| dx'$. Changing variables by $y' = \delta x'$, we obtain

$$\frac{K_\eta}{M_\eta} = \frac{K_1 \eta^{-(n-2)}}{M_1 \eta^{-(n-1)}} = \frac{K_1}{M_1} \eta,$$

where $M_1 = \int_{B_1} |\phi(y')| dy'$ and $K_1 = \int_{B_1} |\nabla_{y'} \phi(y')| dy'$. Thus

$$Q_\gamma[u_{\varepsilon,\eta}] = \frac{\frac{K_1}{M_1} \eta \int_\varepsilon^\delta x_n^{1-s} dx_n + 2\delta^{1-s}}{\int_\varepsilon^\delta x_n^{-1} X^\gamma(x_n/R) dx_n}.$$

Now we select $\eta = \varepsilon^{s-2+\epsilon}$ for some fixed $\epsilon > 0$. We deduce

$$Q_1[u_{\varepsilon,\eta}] = \frac{\frac{K_1}{M_1} \varepsilon^{s-2+\epsilon} \int_\varepsilon^\delta x_n^{1-s} dx_n + 2\delta^{1-s}}{\log(\frac{X(\delta/R)}{X(\varepsilon/R)})}.$$

It follows that $Q_1[u_{\varepsilon,\eta}] \to 0$, as $\varepsilon \downarrow 0$. Thus, for $\Omega = S_R$ inequality (2.10) does not hold when $\gamma = 1$ and the exponent 1 on the distance function in the remainder term in (2.10) cannot be increased.

# 5 Proof of Theorem C

In this section, we assume $\Omega$ is a ball of radius $R$. Without loss of generality, we assume it is centered at the origin, and denote it by $B_R$. The distance function to the boundary is then $d(x) = R - r$, where $r := |x|$. Moreover,

$$-\Delta d(x) = \frac{n-1}{R - d(x)}, \quad x \in B_R \setminus \{0\}. \tag{5.1}$$

This section is devoted to the proof of the following fact

**Theorem 5.1.** (1) *For all* $u \in C_c^\infty(B_R)$, $s \geq 2$ *and* $\gamma > 1$, *it holds that*

$$\int_{B_R} \frac{|\nabla u|}{d^{s-1}} dx \geq (s-1) \int_{B_R} \frac{|u|}{d^s} dx + \sum_{k=1}^{[s]-1} \frac{n-1}{R^k} \int_{B_R} \frac{|u|}{d^{s-k}} dx + \frac{C}{R^{s-1}} \int_{B_R} \frac{|u|}{d} X^\gamma\left(\frac{d}{R}\right) dx, \tag{5.2}$$

*where* $C \geq \gamma - 1$. *The exponents* $s - k$; $k = 1, 2, ..., [s] - 1$, *on the distance function as well as the constants* $(n-1)/R^k$; $k = 1, 2, ..., [s] - 1$, *in the summation terms are optimal. If* $\gamma = 1$ *the above inequality fails in the sense of* (5.5).

(2) *For all* $u \in C_c^\infty(B_R)$, $1 \leq s < 2$ *and* $\gamma > 1$, *it holds that*

$$\int_{B_R} \frac{|\nabla u|}{d^{s-1}} dx \geq (s-1) \int_{B_R} \frac{|u|}{d^s} dx + \frac{C}{R^{s-1}} \int_{B_R} \frac{|u|}{d} X^\gamma\left(\frac{d}{R}\right) dx, \tag{5.3}$$

*where* $C \geq \gamma - 1$. *If* $\gamma = 1$ *the above inequality fails in the sense of* (5.5).

**Remark 5.2.** The optimality of the exponents and the constants stated in the above Theorem is meant in the following sense: for any $s \geq 1$ set

$$I_0[u] := \int_{B_R} \frac{|\nabla u|}{d^{s-1}} dx - (s-1) \int_{B_R} \frac{|u|}{d^s} dx,$$

and also for any $s \geq 2$ set

$$I_m[u] := I_0[u] - \sum_{k=1}^{m} \frac{n-1}{R^k} \int_{B_R} \frac{|u|}{d^{s-k}} dx, \quad m = 1, ..., [s] - 1.$$

Then, for any $s \geq 2$

$$\inf_{u \in C_c^\infty(B_R) \setminus \{0\}} \frac{I_m[u]}{\int_{B_R} \frac{|u|}{d^\beta} dx} = \begin{cases} (n-1)/R^{m+1}, & \text{if } \beta = s - m - 1 \\ 0, & \text{if } \beta > s - m - 1, \end{cases} \tag{5.4}$$

for all $m \in \{0, ..., [s] - 2\}$. Further, for any $s \geq 1$

$$\inf_{u \in C_c^\infty(B_R) \setminus \{0\}} \frac{I_{[s]-1}[u]}{\int_{B_R} \frac{|u|}{d} X(d/R) dx} = 0. \tag{5.5}$$

*Proof.* Inequality (5.3) is evident by Theorem 2.11. Let $s \geq 2$ and $\gamma > 1$. Since inequality (5.2) is scale invariant it suffices to prove it for $R = 1$. Testing (2.1) with

$$T(x) = -(d(x))^{1-s}[1 - (d(x))^{s-1} X^{\gamma-1}(d(x))]\nabla d(x), \quad x \in B_1 \setminus \{0\}.$$

we arrive at

$$\int_{B_1} \mathrm{div}(T)|u|\mathrm{d}x = (s-1)\int_{B_1} \frac{|u|}{d^s}\mathrm{d}x + \int_{B_1} \frac{|u|}{d^{s-1}}(1 - d^{s-1}X^{\gamma-1}(d))(-\Delta d)\mathrm{d}x$$
$$+ (\gamma - 1)\int_{B_1} \frac{|u|}{d}X^{\gamma}(d)\mathrm{d}x.$$

Thus, using (5.1) for $R = 1$, we obtain

$$\int_{B_1} \mathrm{div}(T)|u|\mathrm{d}x = (s-1)\int_{B_1} \frac{|u|}{d^s}\mathrm{d}x + (n-1)\int_{B_1} \frac{|u|}{d^{s-1}}\frac{1 - d^{s-1}X^{\gamma-1}(d)}{1-d}\mathrm{d}x$$
$$+ (\gamma - 1)\int_{B_1} \frac{|u|}{d}X^{\gamma}(d)\mathrm{d}x. \tag{5.6}$$

Since $s \geq 2$, we take into account in (5.6) the fact that

$$\frac{1 - d^{s-1}X^{\gamma-1}(d)}{1-d} \geq \frac{1 - d^{s-1}}{1-d} \geq \frac{1 - d^{[s]-1}}{1-d} = \sum_{k=1}^{[s]-1} d^{k-1}, \quad x \in B_1 \setminus \{0\},$$

and finally arrive at

$$I_0[u] \geq (s-1)\int_{B_1} \frac{|u|}{d^s}\mathrm{d}x + (n-1)\sum_{k=1}^{[s]-1}\int_{B_1} \frac{|u|}{d^{s-k}}\mathrm{d}x + (\gamma-1)\int_{B_1} \frac{|u|}{d}X^{\gamma}(d)\mathrm{d}x,$$

which is (5.2) for $R = 1$.

We next prove (5.4). Suppose first that $2 \leq s < 3$. In this case all we have to prove is that

$$\inf_{u \in C_0^1(B_1) \setminus \{0\}} \frac{I_0[u]}{\int_{B_1} \frac{|u|}{d^\beta}\mathrm{d}x} = \begin{cases} n-1, & \text{if } \beta = s-1 \\ 0, & \text{if } \beta > s-1. \end{cases} \tag{5.7}$$

To this end, we pick $u_\delta(x) = \chi_{B_{1-\delta}}(x)$, where $x \in B_1$ and $0 < \delta < 1$. This function is in $BV(B_1)$ and we can take a $C_c^\infty$ approximation of it, so that the calculations bellow to hold in the limit. The distributional gradient of $u_\delta$ is $\nabla u_\delta = -\vec{\nu}_{\partial B_{1-\delta}}\delta_{\partial B_{1-\delta}}$, and the total variation of $\nabla u_\delta$ is $|\nabla u_\delta| = \delta_{\partial B_{1-\delta}}$. Using co-area formula we get

$$\frac{I_0[u_\delta]}{\int_{B_1} \frac{|u_\delta|}{d^\beta}\mathrm{d}x} = \frac{\delta^{1-s}|\partial B_{1-\delta}| - (s-1)\int_0^{1-\delta}(1-r)^{-s}|\partial B_r|\mathrm{d}r}{\int_0^{1-\delta}(1-r)^{-\beta}|\partial B_r|\mathrm{d}r}$$
$$= \frac{\delta^{1-s}(1-\delta)^{n-1} - \int_0^{1-\delta}((1-r)^{1-s})'r^{n-1}\mathrm{d}r}{\int_0^{1-\delta}(1-r)^{-\beta}r^{n-1}\mathrm{d}r}$$
$$= (n-1)\frac{\int_0^{1-\delta}(1-r)^{1-s}r^{n-2}\mathrm{d}r}{\int_0^{1-\delta}(1-r)^{-\beta}r^{n-1}\mathrm{d}r}.$$

Thus

$$\frac{I_0[u_\delta]}{\int_{B_1} \frac{|u_\delta|}{d^\beta} dx} \to \begin{cases} n-1, & \text{if } \beta = s-1 \\ 0, & \text{if } \beta > s-1 \end{cases}, \quad \text{as } \delta \downarrow 0.$$

Assume next that $3 \leq s < 4$. This time, besides (5.7) we have to prove that

$$\inf_{u \in C_c^\infty(B_1) \setminus \{0\}} \frac{I_1[u]}{\int_{B_1} \frac{|u|}{d^\beta} dx} = \begin{cases} n-1, & \text{if } \beta = s-2 \\ 0, & \text{if } \beta > s-2. \end{cases}$$

Picking the same $u_\delta$ as before and performing the same integration by parts in the second term of the numerator, we conclude

$$\begin{aligned}
\frac{I_1[u_\delta]}{\int_{B_1} \frac{|u_\delta|}{d^\beta} dx} &= \frac{(n-1)\int_0^{1-\delta}(1-r)^{1-s}r^{n-2}dr - (n-1)\int_0^{1-\delta}(1-r)^{1-s}r^{n-1}dr}{\int_0^{1-\delta}(1-r)^{-\beta}r^{n-1}dr} \\
&= (n-1)\frac{\int_0^{1-\delta}(1-r)^{2-s}r^{n-2}dr}{\int_0^{1-\delta}(1-r)^{-\beta}r^{n-1}dr}.
\end{aligned}$$

Thus

$$\frac{I_1[u_\delta]}{\int_{B_1} \frac{|u_\delta|}{d^\beta} dx} \to \begin{cases} n-1, & \text{if } \beta = s-2 \\ 0, & \text{if } \beta > s-2, \end{cases}, \quad \text{as } \delta \downarrow 0.$$

We continue in the same fashion for $4 \leq s < 5$, then $5 \leq s < 6$ and so on.

Next we prove (5.5). We pick $u_\delta$ as before and perform the same integration by parts, to get

$$\begin{aligned}
\frac{I_{[s]-1}[u_\delta]}{\int_{B_1} \frac{|u_\delta|}{d} X(d) dx} &= \frac{(n-1)\int_0^{1-\delta}(1-r)^{1-s}r^{n-2}dr - (n-1)\sum_{k=1}^{[s]-1}\int_0^{1-\delta}(1-r)^{k-s}r^{n-1}dr}{\int_0^{1-\delta}(1-r)^{-1}r^{n-1}X(1-r)dr} \\
&= (n-1)\frac{\int_0^{1-\delta}(1-r)^{[s]-s}r^{n-2}dr}{\int_1^{1-\log\delta} t^{-1}(1-e^{1-t})^{n-1}dt} \\
&=: (n-1)\frac{N_\delta}{D_\delta}.
\end{aligned}$$

Since $[s] - s > -1$, we have $N_\delta = O_\delta(1)$ as $\delta \downarrow 0$. Also, $D_\delta \geq \int_1^{1-\log\delta} t^{-1}dt + O_\delta(1) \to \infty$, as $\delta \downarrow 0$. ∎

# 6 From $L^1$ to $L^p$ weighted Hardy inequalities

In this section we discuss how far our results can go in the $L^p$ setting. We start with the $L^p$ analog of Lemma 2.2.

**Lemma 6.1.** *Let $\Omega \subsetneq \mathbb{R}^n$ be open. For all $u \in C_c^\infty(\Omega)$, all $s > 1$, $p \geq 1$, it holds that*

$$\int_\Omega \frac{|\nabla u|^p}{d^{s-p}} dx \geq \left(\frac{s-1}{p}\right)^p \int_\Omega \frac{|u|^p}{d^s} dx + \left(\frac{s-1}{p}\right)^{p-1} \int_\Omega \frac{|u|^p}{d^{s-1}}(-\Delta d) dx, \quad (6.1)$$

*where $-\Delta d$ is meant in the distributional sense.*

*Proof.* We substitute $u$ by $|u|^p$ with $p > 1$ in (2.3), to arrive at

$$\frac{p}{s-1} \int_\Omega \frac{|\nabla u||u|^{p-1}}{d^{s-1}} dx \geq \int_\Omega \frac{|u|^p}{d^s} dx + \frac{1}{s-1} \int_\Omega \frac{|u|^p}{d^{s-1}}(-\Delta d) dx. \quad (6.2)$$

The left hand side in (6.2) can be written as follows

$$\frac{p}{s-1} \int_\Omega \frac{|\nabla u||u|^{p-1}}{d^{s-1}} dx = \int_\Omega \left\{\frac{p}{s-1} \frac{|\nabla u|}{d^{s/p-1}}\right\} \left\{\frac{|u|^{p-1}}{d^{s-s/p}}\right\} dx$$

$$\leq \frac{1}{p}\left(\frac{p}{s-1}\right)^p \int_\Omega \frac{|\nabla u|^p}{d^{s-p}} dx + \frac{p-1}{p} \int_\Omega \frac{|u|^p}{d^s} dx,$$

by Young's inequality. Thus (6.2) becomes

$$\frac{1}{p}\left(\frac{p}{s-1}\right)^p \int_\Omega \frac{|\nabla u|^p}{d^{s-p}} dx \geq \frac{1}{p} \int_\Omega \frac{|u|^p}{d^s} dx + \frac{1}{s-1} \int_\Omega \frac{|u|^p}{d^{s-1}}(-\Delta d) dx.$$

Rearranging the constants we arrive at the inequality we sought for. ∎

**Remark 6.2.** (**I**) If $\Omega$ satisfies condition (C), we may cancel the last term to obtain

$$\int_\Omega \frac{|\nabla u|^p}{d^{s-p}} dx \geq \left(\frac{s-1}{p}\right)^p \int_\Omega \frac{|u|^p}{d^s} dx, \quad (s > 1, \ p \geq 1).$$

The constant is optimal, as can be seen by arguing as in the proof of Theorem 2.11, with the choice $u_\varepsilon(x) = (d(x))^{(s-1)/p+\varepsilon}\phi(x) \in W_0^{1,p}(\Omega; d^{-(s-p)}); \varepsilon > 0$, and using the elementary inequality $|a+b|^p \leq |a|^p + c_p(|a|^{p-1}|b| + |b|^p); a, b \in \mathbb{R}^n$ and $p > 1$, in the numerator. Under the stronger assumption that $\Omega$ is convex this result was given in [Avkh] by approximation with bounded convex polytopes, while in the non weighted case, i.e., $s = p$, and $\Omega$ convex, it was first given in [MMP]. See also [MS] for the two dimensional non weighted case.

(**II**) Inequality (6.1) in the non-weighted case is proved in [BFT1]-Lemma 3.3(ii) where in addition an extra term appears on the right hand side. This allowed the authors to obtain even more singular potentials for domains having finite inner radius and satisfying property (C). In particular, the optimal homogeneous improvement was obtained (see [3-4]-Theorem A).

(**III**) If $\Omega$ is bounded, the second constant appearing on the right hand side in (6.1) is optimal. To see this, we choose $u_\varepsilon(x) = (d(x))^{(s-1)/p+\varepsilon} \in W_0^{1,p}(\Omega; d^{-(s-p)})$, $\varepsilon > 0$, and after simple computations, involving an integration by parts in the denominator, we conclude

$$\frac{\int_\Omega \frac{|\nabla u_\varepsilon|^p}{d^{s-p}} dx - \left(\frac{s-1}{p}\right)^p \int_\Omega \frac{|u_\varepsilon|^p}{d^s} dx}{\int_\Omega \frac{|u_\varepsilon|^p}{d^{s-1}}(-\Delta d) dx} = \frac{\left(\frac{s-1}{p}+\varepsilon\right)^p - \left(\frac{s-1}{p}\right)^p}{\varepsilon p}$$

$$\to \left(\frac{s-1}{p}\right)^{p-1}, \quad \text{as } \varepsilon \downarrow 0.$$

(**IV**) By (I) and (III), if $\Omega$ is a bounded set such that condition (C) holds, then all constants appearing in (6.1) are optimal.

(**V**) Assume finally that $\Omega$ is a domain with $\partial\Omega \in \mathcal{C}^2$ and such that it satisfies a uniform interior sphere condition. By Theorem 3.4 and (6.1), we get the $L^p$ analog of Theorem 4.2

$$\int_\Omega \frac{|\nabla u|^p}{d^{s-p}}\mathrm{d}x \geq \Big(\frac{s-1}{p}\Big)^p \int_\Omega \frac{|u|^p}{d^s}\mathrm{d}x + (n-1)\underline{\mathcal{H}}\Big(\frac{s-1}{p}\Big)^{p-1} \int_\Omega \frac{|u|^p}{d^{s-1}}\mathrm{d}x.$$

**Acknowledgements** I would like to thank my PhD supervisor, Prof. Stathis Filippas, for his essential help in this work.

# References

[AltC] ALTER, F., CASELLES, V.: *Uniqueness of the Cheeger set of a convex body.* Nonlinear Anal. **70** (1), 32-44 (2009).

[AmbM] AMBROSIO, L., MANTEGAZZA, C.: *Curvature and distance function from a manifold.* J. Geom. Anal. **8** (5), 723 - 748 (1997).

[Avkh] AVKHADIEV, F. G.: *Hardy type inequalities in higher dimensions with explicit estimate of constants.* Lobachevskii J. Math. **21**, 3-31 (2006).

[BFT1] BARBATIS, G., FILIPPAS. S., TERTIKAS, A.: *A unified approach to improved $L^p$ Hardy inequalities with best constants.* Trans. Amer. Math. Soc. **356** (6), 2169-2196 (2003).

[BFT2] BARBATIS, G., FILIPPAS. S., TERTIKAS, A.: *Series expansion for $L^p$ Hardy inequalities.* Indiana Univ. Math. J. **52** (1), 171-190 (2003).

[BFT3] BARBATIS, G., FILIPPAS. S., TERTIKAS, A.: *Refined geometric $L^p$ Hardy inequalities.* Commun. Contemp. Math. **5** (6), 869-881 (2003).

[BrM] BREZIS, H., MARCUS, M.: *Hardy's inequalities revisited.* Ann. Sc. Norm. Super. Pisa Cl. Sci. **25** (1-2), 217–237 (1997).

[BZ] BURAGO, YU. D., ZALGALLER, V. A.: **Geometric Inequalities**. Grundlehren der mathematischen Wissenschaften **285**, Springer (1988).

[CC] CANNARSA, P., CARDALIAGUET, P.: *Representation of equilibrium solutions to the table problem for growing sandpiles.* J. Eur. Math. Soc. **6**, 1-30 (2004).

[CS] CANNARSA, P., SINESTRARI, C.: **Semiconcave Functions, Hamilton-Jacobi Equations, and Optimal Control**. Progress in Nonlinear Differential Equations and Their Applications **58**, Birkhäuser Boston, Inc. (2004).

[CrM] CRASTA, G., MALUSA, A.: *The distance function from the boundary in a Minkowski space.* Trans. Amer. Math. Soc. **359** (12), 5725-5759 (2007).

[D] DAVIES, E. B.: *Some norm bounds and quadratic form inequalities for Schrodinger operators, II.* J. Operator Theory **12** (1), 177-196 (1984).

[EvG] EVANS, L. C., GARIEPY, R. F.: **Measure theory & fine properties of functions**. Studies in Advanced Mathematics, CRC Press (1992).

[F] FEDERER, H.: *Curvature measures.* Trans. Amer. Math. Soc. **93** (3), 418-491 (1959).

[FMT1] FILIPPAS, S., MAZ'YA, V. G., TERTIKAS, A.: *A sharp Hardy Sobolev inequality.* C. R. Math. Acad. Sci. Paris **339**, 483-486 (2004).

[FMT2] FILIPPAS, S., MAZ'YA, V. G., TERTIKAS, A.: *On a question of Brezis and Marcus.* Calc. Var. Partial Differential Equations **25** (4), 491-501 (2006).

[FMT3] FILIPPAS, S., MAZ'YA, V. G., TERTIKAS, A.: *Critical Hardy-Sobolev inequalities.* J. Math. Pures Appl. **87** (1), 37–56 (2007).

[FTT] FILIPPAS, S., TERTIKAS, A., TIDBLOM, J.: *On the structure of Hardy-Sobolev-Maz'ya inequalities.* J. Eur. Math. Soc. **11** (6), 1165-1185 (2009).

[FrK] FRIDMAN V., KAWOHL, B.: *Isoperimetric estimates for the first eigenvalue of the p-Laplace operator and the Cheeger constant.* Comment. Math. Univ. Carolinae **44** (4), 659-667 (2003).

[Fu] FU, J. H. G.: *Tubular neighborhoods in Euclidean spaces.* Duke Math. J. **52** (4), 1025-1046 (1985).

[GTr] GILBARG D., TRUDINGER, N. S.: **Elliptic Partial Differential Equations of second order** (2nd edition). Grundlehren der mathematischen Wissenschaften **224**, Springer (1983).

[G] GIORGIERI, E.: A boundary value problem for a PDE model in mass transfer theory: Representation of solutions and regularity results. PhD Thesis, Universitá di Roma Tor Vergata, Rome, Italy, 2004.

http://art.torvergata.it/bitstream/2108/211/1/tesi-07-12-04.pdf

[GN] GOLDMAN, M., NOVAGA, M.: *Volume-constrained minimizers for the prescribed curvature problem in periodic media.* `arXiv:1103.5161v4`.

[Gr] GROMOV, M.: *Sign and geometric meaning of curvature.* Rend. Semin. Mat. Fis. Milano **61**, 9-123 (1991).

www.ihes.fr/˜gromov/PDF/1[77].pdf

[HLP] HARDY, G., LITTLEWOOD, J. E., PÓLYA, G.: **Inequalities**. Cambridge University Press (1934).

[LL] LEWIS, R. T., LI, J.: *A geometric characterization of a sharp Hardy inequality.* `arXiv:1103.5429`.

[LLL] LEWIS, R. T., LI, J., LI, Y.-Y.: *A geometric characterization of a sharp Hardy inequality.* J. Funct. Anal. **262** (7), 3159-3185 (2012).

[LN] LI, Y.Y., NIRENBERG, L.: *Regularity of the distance function to the boundary.* Rend. Accad. Naz. Sci. XL Mem. Mat. Appl. **29**, 257-264 (2005).

[MMP] MARCUS, M., MIZEL, V. J., PINCHOVER, Y.: *On the best constant for Hardy's inequality in $\mathbb{R}^n$.* Trans. Amer. Math. Soc. **350** (8), 3237-3255 (1998).

[MS] MATSKEWICH, T., SOBOLEVSKII, P.E.: *The best possible constant in generalized Hardy's inequality for convex domain in $\mathbb{R}^n$.* Nonlinear Anal. **28** (9), 1601-1610 (1997).

[Mz] MAZ'YA, V. G.: **Sobolev Spaces**. Translated from the Russian by T. O. Shaposhnikova. Springer Series in Soviet Mathematics, Springer (1985).

[Mnn] MENNUCCI, A. C. G.: *Regularity and variationality of solutions to Hamilton-Jacobi equations Part I: regularity.* ESAIM Control Optim. Calc. Var. **10**, 426-451 (2004) & *Errata.* ESAIM Control Optim. Calc. Var. **13** (2), 413-417 (2007).

[S] SANTALÓ, L. A.: **Integral geometry & Geometric Probability**. Encyclopedia of Mathematics and its Applications **1**, Addison-Wesley Publishing Co. (1979).

[StrZ] STREDULINSKY, E., ZIEMER, W. P.: *Area minimizing sets subject to a volume constraint in a convex set.* J. Geom. Anal. **7** (4), 653 - 677 (1997).